\RequirePackage{fix-cm}
\documentclass[twocolumn]{svjour3}          % twocolumn
\smartqed  % flush right qed marks, e.g. at end of proof
\usepackage{graphicx}
\usepackage{amsfonts}
\usepackage{amsmath}
\usepackage{subfigure}
\usepackage{algorithmic}
\usepackage{algorithm}
\usepackage{slashbox}
\usepackage{epstopdf}
\usepackage{color}
\begin{document}

\title{Well Control Optimization using Derivative-Free Algorithms and a Multiscale Approach}


\titlerunning{Multiscale Approach}        % if too long for running head

\author{Xiang Wang         \and
        Ronald D. Haynes   \and
        Qihong Feng   %etc.
}

%\authorrunning{Multiscale Approach} % if too long for running head

\institute{X. WANG \at
              Department of Petroleum Engineering, China University of Petroleum (Huadong), Qingdao, Shandong, China\\
              Department of Mathematics \& Statistics, Memorial University, St. John's, NL, Canada\\
              Tel.: +86.132.8085.6026\\
              \email{xiangwangdr@gmail.com}           %  \\
           \and
           R.D. Haynes \at
              Department of Mathematics \& Statistics, Memorial University, St. John's, NL, Canada \\
              Tel.: +1.709.864.8825\\
              Fax: +1.709.864.3010\\
              \email{rhaynes@mun.ca}
            \and
          Q. Feng \at
              Department of Petroleum Engineering, China University of Petroleum (Huadong), Qingdao, Shandong, China \\
              \email{fengqihong@126.com}
}

\date{Received: date / Accepted: date}

\maketitle

\begin{abstract}
Smart well technologies, which allow remote control of well and production processes, make the problem of determining optimal control strategies a timely endeavour. In this paper, we use numerical optimization algorithms and a multiscale approach in order to find an optimal well management strategy over the life of the reservoir.  Optimality is measured in terms of the values of the net present value objective function. The large number of well rates for each control step make the optimization problem more difficult and  at a high risk of achieving a suboptimal solution. Moreover, the optimal number of adjustments is not known a priori. Adjusting well controls too frequently will increase unnecessary well management and operation cost, and an excessively low number of control adjustments may not be enough to obtain a good yield. We investigate three derivative-free optimization algorithms, chosen for  their robust and parallel nature, to determine  optimal well control strategies. The algorithms chosen  include generalized pattern search (GPS), particle swarm optimization (PSO) and covariance matrix adaptation evolution strategy (CMA-ES). These three algorithms encompass the breadth of available black--box optimization strategies: deterministic local search, stochastic global search and stochastic local search. In addition, we hybridize the three derivative-free algorithms with a multiscale regularization approach. Starting with a reasonably small number of control steps, the control intervals are subsequently refined during the optimization. Results for experiments studied indicate that CMA-ES performs best among the three algorithms in solving both small and large scale problems. When hybridized with a multiscale regularization approach, the ability to find the optimal solution is further enhanced, with the performance of GPS improving the most.
Topics affecting the performance of the multiscale approach are discussed in this paper, including the effect of control frequency on the well control problem. The parameter settings for GPS, PSO, and CMA-ES, within the multiscale approach are considered.

\end{abstract}

\section{Introduction}
\label{sec:1}

Determining the well production and injection rates is of paramount importance in modern reservoir development. The decision is difficult since the optimal rates depend on the heterogeneity of the rock and liquids, the well placements and other parameters.  Indeed, these properties and input parameters are coupled in a highly nonlinear fashion. Moreover, the optimal production and injection rates are usually not constant throughout the life cycle of reservoir. The oil saturation distribution changes during the well injection and production processes. This will then affect the optimal  production and injection rate for each well. 

Well control planning can be formulated as an optimization problem, using economic or cumulative oil production as the objective function. The well rates or bottom hole pressures at different times are the optimization variables. Many optimization algorithms have been investigated to solve such problems. These algorithms can be broadly placed in two categories: derivative-based algorithms and derivative-free algorithms \cite{isebor_constrained_2009}.

Derivative-based or gradient-based algorithms, take advantage of the gradient information to guide their search. This type of algorithm, commonly used in well control optimization, includes steepest ascent \cite{wang_production_2009}, conjugate gradient \cite{aitokhuehi_real-time_2004}, and sequential quadratic programming  methods \cite{isebor_constrained_2009}. Gradients of the objective function may be calculated by using an adjoint equation. This is an invasive approach, requiring a  detailed knowledge of mathematics inside the reservoir simulator  \cite{isebor_constrained_2009,brouwer_dynamic_2004}. Other ways to approximate the gradients include methods such as finite difference perturbation \cite{aitokhuehi_real-time_2004}, or the simultaneous perturbation stochastic approximation \cite{wang_production_2009}. These algorithms assume a certain degree of smoothness of the objective function with respect to the optimization variables.  Derivative-based algorithms are potentially very quick to converge but sometimes fall into local optimal. 

Derivative-free algorithms can be subdivided into local search methods and global search methods. Local derivative-free algorithms include generalized pattern search (GPS)  \cite{humphries_simultaneous_2013}, mesh adaptive direct search (MADS) \cite{isebor_constrained_2009}, Hooke-Jeeves direct search (HJDS)  \cite{isebor_constrained_2009}, ensemble-based optimization (EnOpt) \cite{oliveira_adaptive_2014}, covariance matrix adaptation evolution strategy (CMA-ES) \cite{bouzarkouna_well_2012,loshchilov_cma-es_2013}, and so on. These methods have strong ability to find accurate optima in a local space, but may face with some difficulties in finding global optima, especially when a good initial guess is not available. Global derivative-free algorithms search through the entire space and provide techniques to avoid being trapped in local optima. Examples of global search algorithms include genetic algorithms (GAs) \cite{almeida_evolutionary_2007}, particle swarm optimization (PSO) \cite{isebor_generalized_2014}, and differential evolution (DE) \cite{storn_differential_1997,Grazieli_de_2015}. Although these algorithms are robust and easy to use, they often require more function evaluations than local search and derivative-based algorithms to converge. However, most of these algorithms parallelize naturally and easily, which make their efficiency satisfactory \cite{ciaurri_derivative-free_2011}. Recently, some hybridization of these techniques such as PSO-MADS\cite{humphries_simultaneous_2013,isebor_generalized_2014,humphries_joint_2014}, multilevel coordinate search (MCS) \cite{huyer_global_1999}, etc., are developed and applied in well placement and/or well control optimization problems. These methods provide global search capabilities in addition to local convergence. The performance of MCS for well placement and control optimization is discussed in \cite{xiang_mcs_2015}.

The optimization algorithms mentioned above can be further classified as either stochastic or deterministic. Stochastic methods use information from the previous iteration and a random element to generate new search points. The random component of the algorithm makes it more likely to avoid local optima, but it may also make the control of solution quality difficult, especially with a limited computational budget. The stochastic algorithms from the above list are MADS, CMA-ES, GA, PSO, and DE. Deterministic methods have no random element. For a given problem, deterministic methods will give the same results for each trial. GPS, HJDS, and MCS are examples of deterministic algorithms. 

All the above mentioned algorithms have been used in well control optimization and/or well placement optimization problems. The performance of the algorithms are problem-dependent.
%One of the researchers' work is test the algorithms performance with various problems and try to find a algorithm that performs stable in most cases.
Some of the algorithms, like GPS (1960s), GA (1960s), and PSO (1990s), have been around for decades, and have been used in petroleum industrial problems for a relatively long time. People has accumulated a great deal of experience through case studies. CMA-ES, developed in 2000s \cite{hansen2001ecj}, was first used in petroleum related optimization problems only in 2012 \cite{bouzarkouna_well_2012}. Though CMA-ES showed great performance in solving well placement optimization problem \cite{bouzarkouna_well_2012}; to date there has been little work to apply it in well control problems.

Although many optimization algorithms have been used, well control optimization is still challenging problem and an active area of research. The number of optimization variables is large in many real--life scenarios. The required number of function evaluations will rise sharply with the increase of the number of variables. A single function evaluation requires one reservoir simulation which is often very demanding in terms of CPU time. The non-convex, non-smooth and multi-modal objective surface further increases the optimization difficulty. 

It is difficult to find a reasonable frequency for well control; an excessively low number of control adjustments may not truly optimize oil recovery. Adjusting each well control too frequently imposes an unrealistic control burden on operations, increasing the total well management cost. Moreover, imposing a high number of control adjustments increases the complexity of the optimization problem so much that there is a high risk of optimization algorithms becoming trapped at local optima and hence missing the optimal startegy \cite{shuai_using_2011}. Multiscale regularization approaches are developed to address these problems. The main idea of the multiscale approach is to start the optimization process with a very coarse control frequency (and thus, with a small number of control variables) and refine successively. The solution at the coarse-scale is used as the initial guess of controls for the next finer scale optimization \cite{shuai_using_2011}. 

Lien et al. \cite{lien_multiscale_2008} used an adaptive multiscale regularization approach with gradient-based algorithms for well control optimization. The number and time of control adjustments are chosen based on heuristically defined refinement indicators, which are calculated by using gradients of the objective function. Shuai et al. \cite{shuai_using_2011} applied ordinary multiscale regularization, also called a successive-splitting multiscale approach, to find appropriate frequencies for well control adjustment. The optimization starts with a coarse number of control adjustments  and subsequently splits each control step into two new ones at every iteration. Two optimization algorithms are considered in their works: ensemble-based optimization (EnOpt) \cite{oliveira_adaptive_2014,shuai_using_2011} and bound optimization by quadratic approximation (BOBYQA) \cite{powell_bobyqa_2009}. More recently, Oliveira et al. \cite{oliveira_adaptive_2014} provide an adaptive hierarchical multiscale approach with EnOpt and an adjoint method for estimation of optimal well controls. The number and lengths of controls are selected adaptively by splitting/merging control steps on the optimization proceeds.
%These papers have combined multiscale approaches with several optimization algorithms. The number of algorithms used in multiscale approaches is still limited, and 
Most of the algorithms are gradient-based. The performance and suitability of multiscale methods combined with derivative-free algorithms for well control optimization still needs further attention and is the focus of this paper.

Given the prevalence of the use of derivative-free optimization algorithms for well control optimization and the need for a multiscale approach for problems with a large number of control variables, we consider the natural marriage of the two philosophies. 
In this paper, we combine a multiscale framework with three typical derivative-free optimization algorithms for the well control problem. We choose a deterministic local search method -- generalized pattern search (GPS), a stochastic local search method -- covariance matrix adaptation evolution strategy (CMA-ES), and a stochastic global search method -- particle swarm optimization (PSO). The performance of each algorithm is analyzed with two reservoir models. Although GPS, PSO and CMA-ES are widely used in petroleum engineering and many other areas, to the best of our knowledge there has been no attempt to combine these methods with a multiscale approach to solve the  well control optimization problem. Furthermore, we give detailed discussions on the following topics:

\begin{itemize}
\item[\textbullet] The effects of control frequency on well control optimization, including the effects on the production and the effects on the control strategy.
\item[\textbullet] The performance of GPS, PSO, and CMA-ES with and without the multiscale framework for well control optimization.
\item[\textbullet] The performance of the approaches as a function of computational budget and the effect of a parallel environment.
\item[\textbullet] The best parameter settings for GPS, PSO, and CMA-ES to maximize their performance within a multiscale framework.
\item[\textbullet] The best configuration for the multiscale framework, including the parameter settings and the choice of hybridized algorithm.
\end{itemize}

This paper is structured as follows. Section 2 describes the well control problem formulation. Section 3 gives brief description of derivative--free optimization methods used. Section 4 gives detailed description of the multiscale approach used and the framework of combining the multiscale approach with the three optimization algorithms. In Section 5, we detail the computational methodology and describe the reservoir models used in this paper. Section 6, we present the results and discussion for the experiments. Finally, in Section 7, we provide a summary and conclusions of this work.

\section{The well control optimization problem}
\label{sec:2}
In this section, we describe the well control optimization problem, including the objective function of interest, the control variables and the imposed constraints. 

The typical objective function associated with a well control problem 
evaluates an economic model and takes into account different costs such as 
the price of oil, the costs of the injection and the production of water. 
Another alternative is to use the cumulative oil production or the barrel of oil equivalent \cite{bouzarkouna_well_2012}.
In this work, the objective function of interest is the net 
present value (NPV) of a time series of cash flows. For the two-phase flow of oil and water, the NPV is defined by
\begin{equation} \label{eq:npv}
\begin{split}
\mathrm{NPV} \left( \mathbf{u} \right) = &
\sum\limits_{n = 1}^{N_t} 
{\Bigg[ 
\frac{{\Delta {t_n}}}{{{{\left( {1 + b} \right)}^{\frac{{{t_n}}}{\tau }}}}}
\bigg( {
{r_{gp}} {\mathbf{q} _{gp}^n (\mathbf{u})} + 
{r_{op}} {\mathbf{q} _{op}^n (\mathbf{u})} 
}
}\\
%\Bigg\{\left. \bigg\{\left.
&- {{c_{wp}} {\mathbf{q}_{wp}^n} (\mathbf{u}) } 
- {{c_{wi}} {\mathbf{q}_{wi}^n} (\mathbf{u}) }  
\bigg) \Bigg],
\end{split}
\end{equation}
% \begin{equation} \label{eq:npv}
% \mathrm{NPV} \left( \mathbf{u} \right) = \sum\limits_{n = 1}^{{N_t}} {\left[ {\frac{{\Delta {t_n}}}{{{{\left( {1 + b} \right)}^{\frac{{{t_n}}}{\tau }}}}} {{\left( {r_{gp}} {\mathbf{q} _{gp}^n (\mathbf{u})} + {{r_{op}} {\mathbf{q} _{op}^n (\mathbf{u})}  - {c_{wp}} {\mathbf{q}_{wp}^n} (\mathbf{u}) } - {{c_{wi}} {\mathbf{q}_{wi}^n}(\mathbf{u}) }  \right)} } } \right]}
% \end{equation}
where $\mathbf{u}$ is set of control variables during the reservoir productive lifetime;
$\mathbf{q}_{gp}^n$, $\mathbf{q}_{op}^n$ and $\mathbf{q}_{wp}^n$, respectively, denote the average gas rate, the average oil rate and the average water rate for the \textit{n}th time step;
$\mathbf{q}_{wi}^n$ is the average water-injection rate for the \textit{n}th time step; 
$r_{gp}$ and $r_{op}$ are the gas and oil revenue; 
$c_{wp}$ is the disposal cost of produced water; 
$c_{wi}$ is the water injection cost; 
$N_t$ is total number of time steps; 
$t_n$ is the time at the end of \textit{n}th time step; 
and $\Delta t_n$ is \textit{n}th time step size.
The quantity $\tau$ provides the appropriate 
normalization for $t_n$, e.g., $\tau = 365$ days.
The quantity $b$ is the fractional discount rate.

The optimization variables $\mathbf{u}$ could contain the well bottom hole pressures or the well liquid rates. In this work, we control wells by specifying the liquid rates. The vector $\mathbf{u}$ is an $N_u$-dimensional column vector, where $N_u$ is the total number of well controls. Assuming each well has the same frequency of control steps, then $N_u=N_t \cdot N_w$, where $N_w$ is the total number of wells and $N_t$ is the total number of time steps.

Well control optimization during the reservoir life cycle can be expressed as the following mathematical problem:

\begin{eqnarray} 
\label{eq:pb_npv}
\mathrm{max} & \quad \mathrm{NPV} \left( \mathbf{u} \right),  \\ 
\label{eq:pb_bd}
\mathrm{subject \ to} & \quad \mathbf{u}_{lb} \leq \mathbf{u} \leq \mathbf{u}_{ub} ,\\
\label{eq:pb_ieq}
& \quad \mathbf{c} \left( \mathbf{u} \right) \leq 0,\\
\label{eq:pb_eq}
& \quad \mathbf{e} \left( \mathbf{u} \right) = 0,
\end{eqnarray}
where $\mathrm{NPV} \left( \mathbf{u} \right)$ is the objective function given by 
equation (\ref{eq:npv}). And, in order, equations (\ref{eq:pb_bd}--\ref{eq:pb_eq}) are the bound, inequality, and equality constraints (if any) imposed on the problem. The quantities $\mathbf{u}_{lb}$ and $\mathbf{u}_{ub}$ are the lower and upper bounds for control variables, where the inequality is understood to apply component--wise. There are several methods which could be used to handle constraints in derivative-free optimization algorithms in the literature \cite{bouzarkouna_well_2012}. In general, infeasible individuals can be rejected, penalized, or repaired. Although all three algorithms used in this paper can be adapted to handle the three types of constraints, we restrict ourselves to  bound constraints only.

\section{Overview of the derivative--free optimization algorithms used}
\label{sec:3}

In this section, we briefly describe the derivative--free optimization algorithms considered in this paper: GPS, PSO, and CMA-ES. These are typical derivative--free, black--box optimization algorithms. Each method has distinct characteristics and all have been applied successfully to solve reservoir development problems as mentioned in the introduction.

\subsection{Generalized Pattern Search}
\label{sec:3_gps}

The generalized pattern search (GPS) is an example of a direct search method. GPS is a deterministic local search algorithm. It does not directly use or require the gradient of the objective function to be specified (or even to exist). Hence GPS can be used on functions that lack smoothness, those that are not continuous or differentiable. GPS can be applied in situations when the objective function is rough and multi-modal and hence the gradients do not guide directions of global ascent \cite{asadollahi_production_2014,isebor_constrained_2009}. GPS is guaranteed to converge to locally optimal solutions and can provide useful approximate solutions for some global problems  \cite{audet_analysis_2002,torczon_convergence_1997,yin_extended_2000}. 

A basic generalized pattern search proceeds as follows. Choose an initial point and evaluate the objective function at that point. Then evaluate the function on a pattern specified by set of directions and a step or mesh size. After the search is complete along all the directions of the pattern, the point with the highest function value becomes the current point for the next iteration. The step size for the next iteration will be multiplied by a factor that is larger than 1 (an expansion factor) if one or more previous moves found a better point, or a factor between 0 and 1 (a contraction factor)  if no better point is found. In one iteration, if all pattern directions are evaluated before choosing a new current point, we call it a complete poll. The algorithm can also provide an incomplete poll procedure instead of the complete one. For an incomplete poll, the algorithm stops searching along the directions as soon as it finds a point whose objective function value is more than that of the current point, and the point it finds becomes the current point at the next iteration. The search stops when, for example, a specified minimum step size may be reached.

\begin{figure}[htbp]
  \centering 
  \subfigure[Maximal positive basis]{ 
    \label{fig:subfig:gps_p2n} %% label for first subfigure 
    \includegraphics[width=0.25\textwidth]{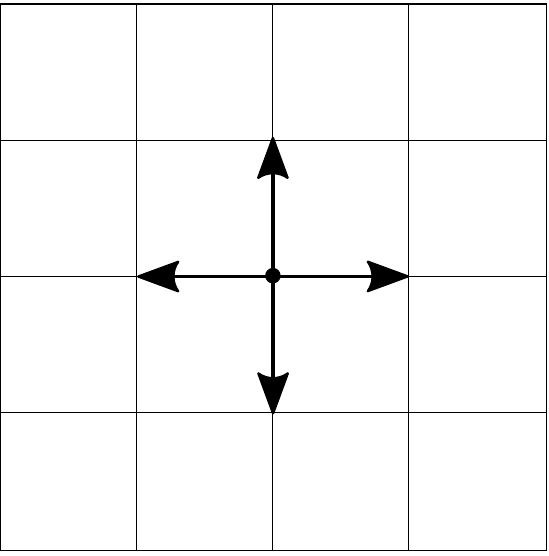}} 
  %\hspace{1in} 
  \subfigure[Minimal positive basis]{ 
    \label{fig:subfig:gps_pn1} %% label for second subfigure 
    \includegraphics[width=0.25\textwidth]{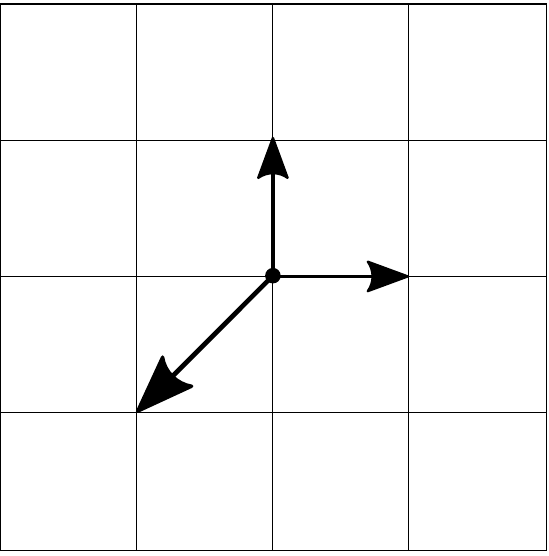}}   
  \caption{Maximal positive basis and minimal positive basis vectors in $\mathbb{R}^2$.} 
  \label{fig:gps_p} %% label for entire figure 
\end{figure}

The choice of the GPS pattern can dramatically affect the performance of the algorithm. From the current point, the pattern determines the search directions for each iteration. The pattern is usually expressed as a set of vectors $\{\mathbf{v}_i \in \mathbb{R}^n:i=1,\cdots,r\}$ which form a positive basis. Every vector in $\mathbb{R}^n$ can be written as a linear combination $a_1\mathbf{v}_1+\cdots+a_r\mathbf{v}_r$ with all coefficients $a_i$ are zero or greater. No vector of the positive basis can be expressed by a positive linear combination of other members of the basis. Using a positive basis is beneficial in GPS because they give small numbers of search directions \cite{abramson_generalized_2003}. Two positive basis are commonly used as search directions in GPS, namely the maximal basis and the minimal basis. For an $n$-dimensional optimization problem, maximal and minimal basis sets have $2n$ and $n+1$ vectors respectively \cite{kolda_generating_2006}. Fig. \ref{fig:gps_p} shows an example of  maximal positive basis (a) and a minimal positive basis (b) in $\mathbb{R}^2$. 

For bound constrained problems, GPS needs a feasible initial point and keeps feasibility of the iterates by rejecting any trial point that is out of the feasible region. For an infeasible trial point, the objective function is not evaluated and set to infinity.

As mentioned above, GPS is guaranteed to converge to locally optimal solutions. Evaluating the objective function at each point in a maximal basis and using complete poll will be quite expensive, however, considering the complexity of well control problems we indeed choose the maximal basis and complete poll for our tests. The expansion factor is set to 2, and the contraction factor is set to 0.5.

\subsection{Particle Swarm Optimization}
\label{sec:3_pso}

Particle swarm optimization (PSO) is a population-based stochastic search method. The PSO search mechanism mimics the social behaviour of biological organisms such as a flock of birds \cite{kennedy_particle_2011}. PSO can search very large space of candidate solutions, and the stochastic element of the movement of the population reduces the chance of getting trapped at an unsatisfactory local optimum, however, PSO does not guarantee an optimal solution is ever found. In spite of this, PSO has been successfully applied for both well placement and production optimization \cite{onwunalu_development_2009,onwunalu_new_2011,isebor_generalized_2014}.

PSO initially chooses a population of candidate solutions (called a swarm of particles). These particles move through the search space in search of function improvement according to a random rule which updates each particle's position. Each particle's position $\mathbf{x}_i^k$ and velocity $\mathbf{v}_i^k$ changes during each iteration. 
For objective function $f:\mathbb{R}^n\rightarrow \mathbb{R}^n$, both $\mathbf{x}_i^k$ and $\mathbf{v}_i^k$ are $n$-dimensional vectors. The objective function value  for each particle, $f(x_i^k)$, is evaluated simply using the position of the particle. Each particle's movement is guided by the best position it has found so far, $\mathbf{p}_i^k$, and the best known position of all particles or particles in some neighborhood, $\mathbf{g}_i^k$. 

Following initialization, the PSO algorithm  \cite{kennedy_particle_2011} updates each particle's position and velocity as:

\begin{equation}\label{eq:pso_x}
\mathbf{x}_i^{k+1}=\mathbf{x}_i^k+\mathbf{v}_i^{k+1},
\end{equation}
and
\begin{equation}\label{eq:pso_v}
\mathbf{v}_i^{k+1}=w \mathbf{v}_i^k + c_1 \mathbf{r}_1^k \otimes \left( \mathbf{p}_i^k - \mathbf{x}_i^k \right) + c_2 \mathbf{r}_2^k \otimes \left( \mathbf{g}_i^k - \mathbf{x}_i^k \right).
\end{equation}

Equation \ref{eq:pso_v} include three parts: $\mathbf{v}_i^k$ represents the tendency to continue moving along the particle's current direction and velocity, $\left( \mathbf{p}_i^k - \mathbf{x}_i^k \right)$ represents the tendency to move to the best position found by the particle itself so far, and $\left( \mathbf{g}_i^k - \mathbf{x}_i^k \right)$ represents the tendency to move to the best position found by all particles in its neighborhood. The quantities $w$, $c_1$ and $c_2$ are weighting parameters and $\mathbf{r}_1^k$ and $\mathbf{r}_2^k$ are stochastic $n$-dimensional vectors which are generated from the uniform distribution on $(0, 1)$ during each iteration. The operator $\otimes$ indicates  component--wise multiplication.   The random element helps ensure that PSO avoids premature convergence to a local minimum by facilitating sufficient global exploration of the search space \cite{kennedy_particle_2011,vaz_particle_2007}.

Here we use an absorbing strategies to handle bound constraints for PSO. The invalid particles are set to the nearest boundary position. The respective velocity components are set to zero \cite{helwig_particle_2009,clerc_particle_2010}.

Unlike local search methods, such as GPS, the PSO algorithm may avoid local optima with its stochastic and global search capability. PSO has been used in petroleum and other fields with excellent effect. However, no number of function evaluations can guarantee convergence to the global optima. The selection of the algorithmic parameters has a considerable affect on the performance of the algorithm \cite{clerc_stagnation_2006,perez_particle_2007}. For the parameters $w$, $c_1$ and $c_2$, Perez et al. \cite{perez_particle_2007} demonstrated that the particle swarm is only stable if the following two conditions are satisfied:
\begin{equation}
0< (c_1+c_2)<4,
\end{equation}
and
\begin{equation}
\frac{1}{2}(c_1+c_2)-1<w<1.
\end{equation}

Following the above principles, our implementation of PSO uses weighting parameters of $w = 0.9$, $c_1 =0.5$, and $c_2 = 1.25$. We use a global best neighbourhood topology, meaning that every particle communicates with every other particle in the swarm, and thus $\mathbf{g}_i^k$ can be replaced by a single vector $\mathbf{g}^k$, representing the best solution found so far. 
Evaluating the objective function for all members of the swarm is an embarrassingly parallel operation. 

\subsection{Covariance Matrix Adaptation Evolution Strategy}
\label{sec:3_cmaes}

Covariance Matrix Adaptation Evolution Strategy (CMA-ES) is a population-based stochastic optimization algorithm. Unlike GA, PSO, and other classical population-based stochastic search algorithms, candidate solutions  of CMA-ES are sampled from a probability distribution which is updated iteratively. This algorithm performs better on the benchmark multimodal functions than all other similar classes of learning algorithms  \cite{willjuice_iruthayarajan_covariance_2010}. CMA-ES also showed its potential in well placement and control optimizations  \cite{bouzarkouna_well_2012,forouzanfar_covariance_2015}.

CMA-ES samples a population of $\lambda$ candidate solutions at iteration $k$ according to:
\begin{equation}\label{eq:cmaes_x}
{\mathbf{x}}_i^{k} = {\mathcal{N}}\left( {{{\mathbf{m}}^{k}}, ({\sigma ^{k}})^2 {{\mathbf{C}}^{k}}} \right) ,\quad {\rm{for}} \ i = 1, \cdots ,\lambda, 
\end{equation}
where $\mathcal{N}(\cdots,\cdots)$ is a random vector from a multivariate normal distribution.

The mean vector $\mathbf{m}^k$ represents the favorite solution or best estimate of the optimum, and the covariance matrix $\mathbf{C}^k$ is a symmetric positive definite matrix which characterizes the geometric shape of the distribution and defines where new candidate solutions are sampled. The step-size $\sigma$ is used as a global scaling factor for the covariance matrix. It aims at achieving fast convergence and preventing premature convergence. These three parameters are updated as the iteration proceeds.

The $\lambda$ individuals generated by equation (\ref{eq:cmaes_x}) are evaluated and ranked by objective function value. The mean $\mathbf{m}^{k}$ is then updated by taking the weighted mean of the best $\mu$ individuals:
\begin{equation}\label{eq:cmaes_m}
{{\mathbf{m}}^{k+1}} = \sum\limits_{i = 1}^\mu  {{\omega _i}{\mathbf{x}}_{1:\lambda }^{k}}  
\end{equation} 
where ${\mathbf{x}}_{1:\lambda }^{k}$ is the $i^{th}$ best individual. 

The default weights \cite{arnold_optimal_2005} are chosen as 
\begin{equation}
{\omega _i} = \frac{{\ln \left( {\mu  + 1} \right) - \ln \left( i \right)}}{{\mu \ln \left( {\mu  + 1} \right) - \ln \left( {\mu !} \right)}},\quad {\rm{for}} \  i = 1, \cdots ,\mu. 
\end{equation}
Typically $\mu$ is chosen as $\mu=\lfloor\lambda/2\rfloor$, where $\lfloor\ \rfloor$ is the floor function, and $\omega _i$ are strictly positive and normalized weights. 

The covariance matrix ${{\mathbf{C}}^{k}}$ is then updated as 
\begin{equation}
\begin{split}
{{\bf{C}}^{k+1}} = &\left( {1 - {c_{{\mathop{\rm cov}} }}} \right){{\bf{C}}^{k}} + \frac{{{c_{{\mathop{\rm cov}} }}}}{{{\mu _{{\mathop{\rm cov}} }}}}{\bf{p}}_c^{k+1}{{\bf{p}}_c^{(k+1)}}^T \\
&+ {c_{{\mathop{\rm cov}} }}\left( {1 - \frac{1}{{{\mu _{{\mathop{\rm cov}} }}}}} \right) \\
&\times \sum\limits_{i = 1}^\mu  {{\frac{\omega _i}{{\sigma^{(k)}}^2}{\left(\mathbf{x}_{i:\lambda}^{k+1}-\mathbf{m}^{k}\right)}{\left(\mathbf{x}_{i:\lambda}^{k+1}-\mathbf{m}^{k}\right)}^T}},
\end{split}
\end{equation}
where quantity ${\bf{p}}_c^k$ is called the evolution path. It gives a direction where we expect to see good solutions. The evolution path is
given iteratively as

\begin{equation}\label{eq:cmaes_pc}
{\bf{p}}_c^{k+1} = \left( {1 - {c_c}} \right){\bf{p}}_c^{k} + \sqrt {{c_c}\left( {2 - {c_c}} \right){\mu _{{\rm{eff}}}}} \frac{{{{\bf{m}}^{k+1}} - {{\bf{m}}^{k}}}}{{{\sigma ^{k}}}},
\end{equation}
where $c_c$ is a constant in $(0,1]$. The quantity $\mu_{\rm{eff}}=1/\sum\nolimits_{i=1}^\mu\omega_i^2$ denotes the variance effective selection mass. It is a measure characterizing the recombination. From equation (\ref{eq:cmaes_pc}) we can see that the new search direction ${\bf{p}}_c^{k+1}$ is based on the old direction ${\bf{p}}_c^{k}$ and the descent direction $\frac{{{{\bf{m}}^{k+1}} - {{\bf{m}}^{k}}}}{{{\sigma ^{k}}}}$.

The adaptation of the global step size $\sigma^{k+1}$ is given by 
\begin{equation}
{\sigma ^{k+1}} = {\sigma ^{k}}\exp \left[ {\frac{{{c_\sigma }}}{{{d_\sigma }}}\left( {\frac{{\left\| {{\bf{p}}_\sigma ^{k+1}} \right\|}}{{{\rm{E}}\left\| {\mathcal{N}\left( {0,{\bf{I}}} \right)} \right\|}} - 1} \right)} \right]
\end{equation}
which depends on the conjugate evolution path $\mathbf{p}_\sigma^{k+1}$ given by

\begin{equation}
\begin{split}
\mathbf{p}_{\sigma} ^{k+1} = & \left( 1-c_\sigma \right) \mathbf{p} _{\sigma}^{k} \\
&+ \sqrt{c_{\sigma} \left( 2-c_{\sigma} \right) \mu _{\rm{eff}} } {{\sigma^{k}}}^{-1} {\mathbf{C} ^{k}}^{-\frac{1}{2}} \left( \mathbf{m}^{k+1}-\mathbf{m}^{k} \right).
\end{split}
\end{equation}

In combination with covariance matrix adaptation, step-size adaptation enables linear convergence on a wide range  of, even ill-conditioned,  functions \cite{bouzarkouna_well_2012}.

\begin{figure*}[htbp]
\centering 
\subfigure[Iteration 1]{ 
    %\label{fig:subfig:punqs3_tops} %% label for first subfigure 
    \includegraphics[width=0.31\textwidth]{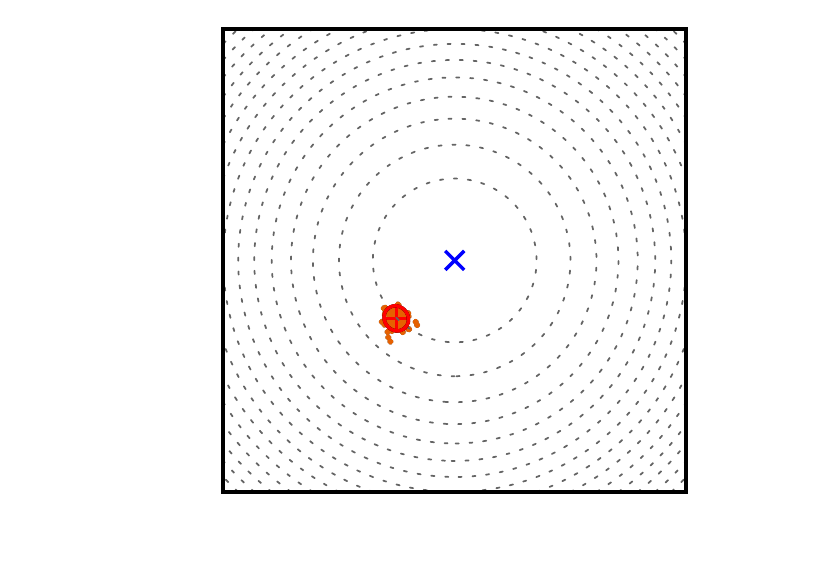}
    }
\subfigure[Iteration 3]{ 
    %\label{fig:subfig:punqs3_perm} %% label for first subfigure 
    \includegraphics[width=0.31\textwidth]{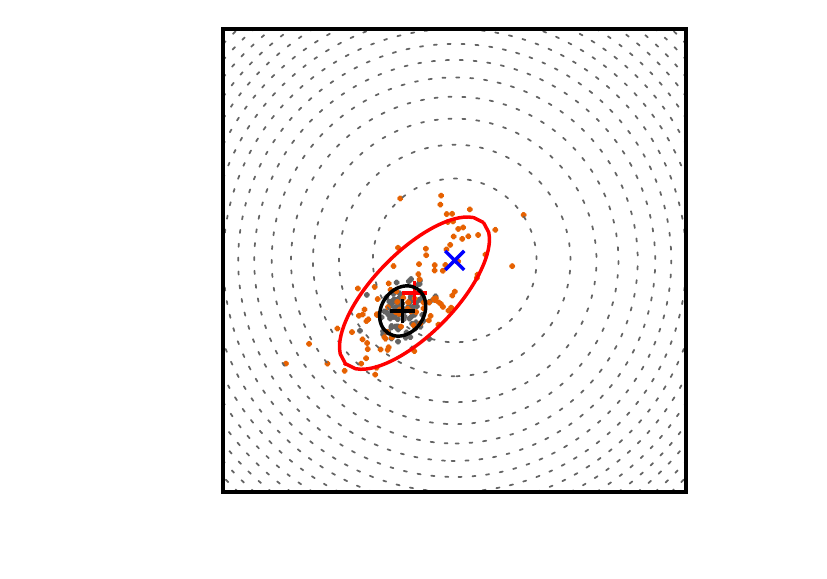}
    }
\subfigure[Iteration 5]{ 
    %\label{fig:subfig:punqs3_soil} %% label for first subfigure 
    \includegraphics[width=0.31\textwidth]{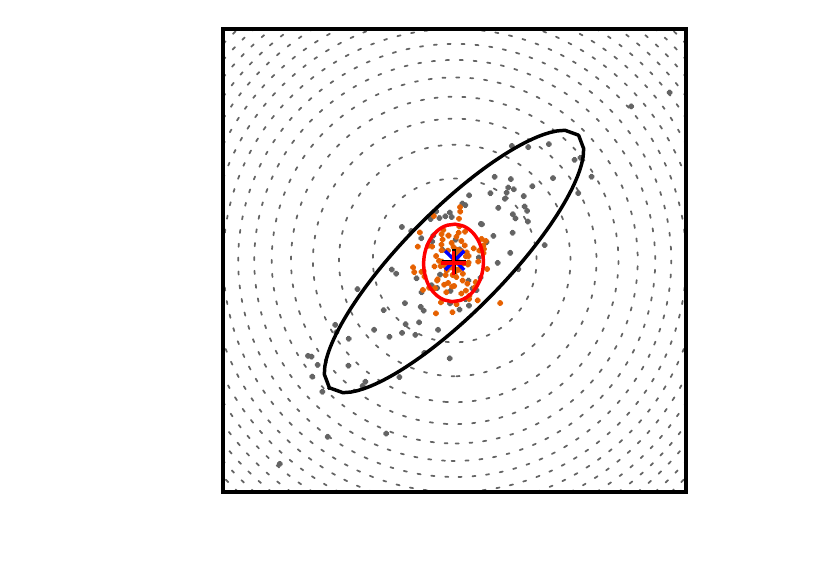}
    }
\caption{An illustration of CMA-ES optimization on a two-dimensional linear function $f(\mathbf{x})=x_1^2+x_2^2$. The dashed lines are the contour lines of the function. The optimal solution is in the upper right corner. The orange and gray dots denote the population distribution for the current and the last iteration, respectively. The black cross ('+') and the black ellipsoid denote the symmetry center and the isodensity line of the distribution for the current iteration, while the gray cross and the gray ellipsoid are for the last iteration.}
\label{fig:cmaes_ex}
\end{figure*}

Fig. \ref{fig:cmaes_ex} is an illustration of an optimization run with CMA-ES on a two-dimensional linear function $f(\mathbf{x})=x_1^2+x_2^2$ with bound constraints $x_1,x_2\in[-800,800]$. The dashed lines are the contour lines of the function. The initial point is $[-200,-200]$. The optimal solution is $[0,0]$ and is marked by blue symbol `x'. The orange and gray dots denote the population distribution for the current and the last iteration, respectively. The red cross (`+') and the red ellipsoid denote the symmetry center $\bf m$ and the isodensity line of the distribution for the current iteration, while the black cross and the black ellipsoid give these quantities for the last iteration. 
The isodensity line defined as $(\mathbf{x}-\mathbf{m})^T\mathbf{C}^{-1}(\mathbf{x}-\mathbf{m})=c$ where the constant $c=3$ \cite{bouzarkouna_well_2012}.
The population is much larger than necessary, but the figure clearly shows how the distribution of the population changes during the optimization.

For bound constrained problems, CMA-ES can simply `repair' an infeasible individual to its nearest feasible solution by using a repair algorithm before the update equations are applied \cite{hansen_method_2009}. This strategy is not recommended because CMA-ES makes implicit assumptions on the distribution of individuals. The distribution can be violated by a repair. Here we choose the penalization strategy from \cite{hansen_cma_2005,bouzarkouna_well_2012}. The infeasible points are repaired by using the repair algorithm and evaluated, and a penalty is added to the function value for every repaired point. The penalty is dependent on the distance to the repaired solution.

Our implementation of CMA-ES uses the parameters given in the work of Hansen et al. \cite{hansen_evaluating_2004}. The parameter values are given in Table \ref{tab:cmaes_para}. The initial values used are $\mathbf{p}_\sigma^{0}=\mathbf{p}_c^{0}=\mathbf{0}$ and $\mathbf{C}^{0}=\mathbf{I}$, while $\mathbf{x}^{(0)}$ and $\sigma^{0}$ are user supplied. 

\begin{table}[htbp]
\centering
\caption{Strategy parameter values used in CMA-ES from \cite{hansen_evaluating_2004}.}
\label{tab:cmaes_para}
\begin{tabular}{cc}
\hline\noalign{\smallskip}
Parameter & Value \\
\noalign{\smallskip}\hline\noalign{\smallskip}
$\lambda$ & $4+\lfloor 3 \ln (n) \rfloor $\\
$\mu$ & $\lfloor \lambda/2 \rfloor$ \\
$c_c $ & $ \frac{4}{n+4} $ \\
$c_\sigma $ & $\frac{\mu_{\rm{eff}}+2}{n+\mu_{\rm{eff}}+3} $ \\
$d_\sigma $ & $1+2\max\left(0,\sqrt{\frac{\mu_{\rm{eff}}-1}{n+1}}-1\right)+c_\sigma $ \\
$\mu_{\rm{cov}} $ & $\mu_{\rm{eff}} $ \\
$c_{\rm{cov}} $ & $\frac{1}{\mu_{\rm{cov}}}\frac{2}{(n+\sqrt{2})^2}+\left(  1-\frac{1}{\mu_{\rm{cov}}}\right)\min\left(1,\frac{2\mu_{\rm{eff}}-1}{(n+2)^2 +\mu_{\rm{eff}}} \right) $ \\
\noalign{\smallskip}\hline
\end{tabular}
\end{table}

\section{A Multiscale framework}
\label{sec:4}

In production optimization, specifying the frequency of needed well control adjustment is a challenge. On one hand, a high frequency adjustment of control parameters imposes unrealistic  burden on operations, leading to an increase in well management costs. In addition, from an optimization perspective a high frequency of control adjustments implies an explosion in the number of control variables, requiring a great amount of computation and time to get an optimal solution. This may be especially true for derivative-free algorithms, which may need many more function evaluations than gradient-based algorithms. Many degrees of freedom also increase the risk of an optimization algorithm being trapped in a local optimum. On the other hand, imposing too few control adjustments  may not truly optimize oil recovery.

Multiscale regularization provides a way to address the complexity of the optimization problem with a large number of control adjustments.  The multiscale approach starts with a coarse number of control steps and successively increases the frequency of control adjustments using the coarse-scale solution as the initial guess for the next finer scale optimization  \cite{lien_multiscale_2008,oliveira_adaptive_2014,shuai_using_2011}. The refinement process is terminated when a specified stopping criteria is satisfied.  For example, a maximum number of control adjustments or a minimum allowable change in the objective function could be imposed. 

To the best of our knowledge, three related multiscale approaches have been investigated for the well control optimization problem. The first approach, first seen in \cite{lien_multiscale_2008}, is referred to as ordinary multiscale or successive-splitting multiscale \cite{shuai_using_2011} . The optimization starts with a coarse number of control adjustments  and subsequently splits each control step into two new ones at every iteration. The second optimization strategy, also proposed by Lien et. al. \cite{lien_multiscale_2008}, uses the magnitude of the components of the gradient of the objective to determine refinement indicators.  The algorithm progressively increases the number of variables using the refinement indicators to choose the  most-efficient partitioning of the current control steps to increase the value of the objective function. The third approach is called the hierarchical multiscale method \cite{oliveira_adaptive_2014,oliveira_hierarchical_2015}. It is similar to the ordinary multiscale approach in \cite{lien_multiscale_2008,shuai_using_2011}, but the algorithm can also merge existing control steps by considering the difference between well controls at two consecutive control steps and the gradient of the objective function with respect to the well controls.
 
The goal of the present work is to explore the feasibility of improving the performance of derivative-free algorithms in solving large scale well control optimization problems by using a multiscale approach. We choose the successive-splitting multiscale approach because the other two methods, the refinement indicator multiscale approach and the hierarchical multiscale approach,  require gradient information of the objective function -- information we do not assume is available.  Furthermore, one recent study compared the sophisticated refinement indicator and hierarchical multiscale approaches with the simpler successive-splitting approach and showed similar performance \cite{oliveira_adaptive_2014}.

Before we introduce our modified approach, we take a look at the original successive-splitting multiscale approach. As Shuai et al. describe in \cite{shuai_using_2011}, the successive-splitting multiscale algorithm generally loops over the following steps:

\begin{enumerate}
\item[1)] INITIALIZATION One control step for each well (initial steps $n_0=1$); The number of unknowns is equal to the number of wells; Initial guesses of control are assigned to each well.
\item[2)] OPTIMIZATION Solving the well control optimization problem using an optimization algorithm.
\item[3)] SPLITTING Split each control step into two steps of equal length (split factor $n_s=2$); This doubles the number of control variables; Use the solution from step 2) as the initial well control; Go to step 2). 
\end{enumerate}

Our experience indicates the efficacy of a multiscale approach depends on two key parameters: the number of control steps for each well at the beginning of the optimization (i.e. the number of initial steps $n_0$) and the multiplicative increase in the number of control steps at every iteration (i.e. the split factor $n_s$). As mentioned, the successive-splitting multiscale approach used in \cite{shuai_using_2011} starts the optimization procedure by finding the optimal control strategy assuming one control step ($n_0=1$). Subsequent optimizations split the number of control steps by a fixed split factor $n_s=2$. We show that this configuration of the two parameters is not always the most efficient configuration. On one hand, the optimal well control strategies with a very coarse parametrization may be dramatically different than with a fine parametrization (or large number of control adjustments). Hence the solution found by a very coarse parametrization is not useful as an initial guess to find the optimal fine parametrization or will require many successive splittings. This observation has to be balanced with the realization and motivation that the problem with a large number of control adjustments is too difficult solve immediately. The split factor is the key to balance the difficulty of optimization problem at each scale and the total number of scales. With a higher split factor, less scales are needed to reach the maximum number of control steps. We will show this is more efficient in some cases.

Based on the above, in addition to coupling the multiscale approach with commonly used derivative free algorithms,  we consider the effect of the choice of the initial number of control number steps $n_0$ and the choice of $n_s$ in the overall efficiency of the multiscale optimization process. We show this added flexibility in our algorithm is useful in some situations. In our modified multiscale approach, we left the choice of the initial number of steps and the choice of the split factor to the user. We start the multiscale algorithm with a reasonably small value of $n_0$ -- the initial number of control steps, and then find the associated optimal controls. 
After maximizing objective function on the basis of the initial control steps, we split each control step into several steps depending on the split factor $n_s$ as 
\begin{equation}
\mathbf{x}_{i+1}(n)=\mathbf{x}_{i*}({\lceil n/n_s \rceil}),\ n=1,2,\cdots,N_w\times n_{s},
\end{equation}
where $\mathbf{x}_{i+1}(n)$ is the $n$th variable in the initial guess for the $(i+1)$th scale; $\mathbf{x}_{i*}({\lceil n/n_s \rceil})$ is the ${\lceil n/n_s \rceil}$th variable in the optimum solution for the $i$th scale, $\lceil\ \rceil$ is the ceiling function; $N_w$ is the total number of variables for the $i$th scale. With this formula, the total number of variables for the $(i+1)$th scale becomes $N_w\times n_{s}$, and every $n_s$ variables for the $(i+1)$th scale use the optimum solution of the $i$th scale.
This process of splitting the control steps and performing a new optimization is continued until the maximum number of control steps is reached.
%For example, if the split factor is $4$ then a control step is split into four steps of equal size for the next scale. This multiples the number of optimization variables  for the next scale by a factor of $4$. 

\begin{figure}[htbp]
\centering
\includegraphics[width=0.45\textwidth]{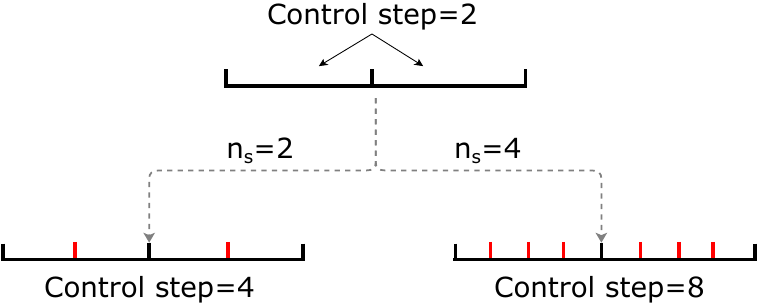}
\caption{Control steps split by the successive-splitting multiscale approach.}
\label{fig:mul_split}
\end{figure}

Fig. \ref{fig:mul_split} gives an illustration of how the successive-splitting multiscale approach splits the control steps to give the next finer scale.  In this figure, we show the resulting number of control steps for two choices of $n_s$ ($n_s=2$ and $n_s=4$) assuming the number of initial steps is $n_0=2$.

\begin{figure}[htbp]
\centering
\includegraphics[width=0.4\textwidth]{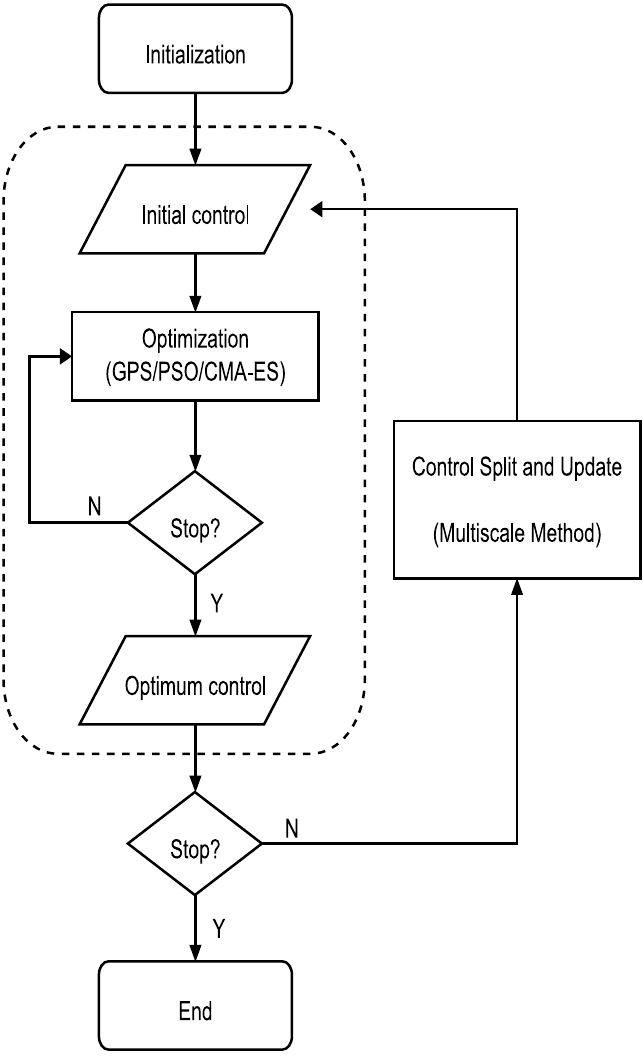}
\caption{Flow chart of well control optimization with multiscale approach and derivative-free algorithms.}
\label{fig:mul_flow}
\end{figure}

Our well control optimization procedure using a derivative-free multiscale approach is described by the following steps (Algorithm \ref{alg1}). A flow chart of the algorithm is given in Fig. \ref{fig:mul_flow}.

\begin{algorithm}                      
\caption{The multiscale approach with derivative-free algorithms}          
\label{alg1}                           
\begin{algorithmic}                    
    \STATE select solver: GPS, PSO \OR CMA-ES
    \STATE set initial control steps for each well $n_0$, and the split factor $n_s$
    \STATE set initial guess $\mathbf{x}_0$
    \STATE iteration $i\leftarrow 0$
	\WHILE{ \NOT (global stopping criteria reached)}
		\WHILE{ \NOT (scale stopping criteria)}
			\STATE solve $ \mathbf{x}_{*} =\mathrm{argmax} \quad \mathrm{NPV} \left( \mathbf{x} \right) $
		\ENDWHILE
		\STATE let $\mathbf{x}_{i*}=\mathbf{x}_{*}$
		\STATE split, set control steps for each well $n_{i+1}\leftarrow n_i\times n_s$
		\STATE update the initial guess, $\mathbf{x}_{i+1}(n)=\mathbf{x}_{i*}({\lceil n/n_s \rceil}),\ n=1,2,\cdots,N_w\times n_{s}$
		\STATE test the scale stopping criteria
    \ENDWHILE
\end{algorithmic}
\end{algorithm}

At early scales of the algorithm, optimizations are performed with less stringent convergence tolerances to find an approximate solution. The tolerances are made smaller as the algorithm proceeds, using the smallest tolerance at the last scale. With this approach, we reduce the computational cost at early scales. The tolerance settings for the experiments can be found in Section \ref{sec:6}.  

\section{Example cases}
\label{sec:5}

In this section, we list all approaches considered, and give a detailed description of the reservoir models used in this paper.

\subsection{Optimization Approaches}
\label{sec:5_1}

Approaches considered in this paper include the three original optimization algorithms, GPS, PSO, and CMA-ES as described in Section \ref{sec:3} and three hybrid approaches that combine the original algorithms with our modified multiscale method described in Section \ref{sec:4}. The hybrid multiscale approaches are labeled as M-GPS, M-PSO, and M-CMA-ES. 

To investigate the effect of $n_0$ and $n_s$, we test four different configurations for each hybrid approach. We use the Roman numerals I, II, III, and IV to represent the four configurations. The configurations used are:  
\begin{itemize}
\item[\textbullet] Configuration I, ---the initial number of control steps for each well is $n_0=1$ and the split factor is $n_s=2$. With this configuration, the multiscale method is the same as the successive-splitting multiscale method from \cite{shuai_using_2011}.
\item[\textbullet] Configuration II, ---the initial number of control steps for each well is $n_0=2$ and the split factor is $n_s=2$. 
\item[\textbullet] Configuration III, ---the initial number of control steps for each well is $n_0=2$ and the split factor is $n_s=4$. 
\item[\textbullet] Configuration IV, ---the initial number of control steps for each well is $n_0=1$ and the split factor is $n_s=4$. 
\end{itemize}

\begin{figure}[htbp]
\centering
\includegraphics[width=0.45\textwidth]{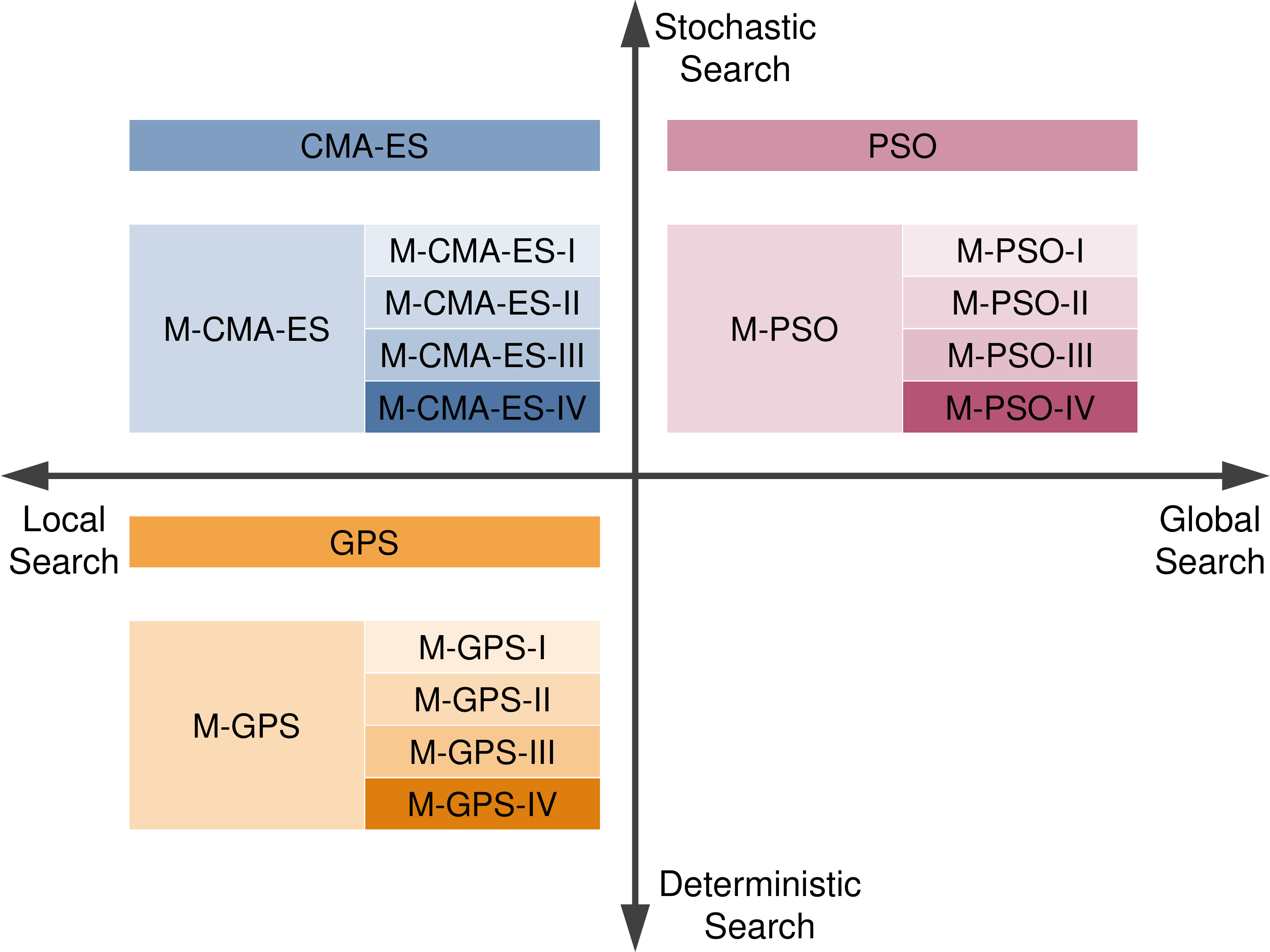}
\caption{An overview of all approaches we considered in our experiments. The approaches fall into different quadrants according to their search features.}
\label{fig:mul_app}
\end{figure}

Fig. \ref{fig:mul_app} provides an overview of all approaches considered in our experiments. The approaches fall into different quadrants according to their search features. 

\subsection{Model description}
\label{sec:5_2}

Two reservoir models are considered in this paper. The first one is a simple 2-D reservoir model. This model is used to analyze the performance of the approaches mentioned in Section \ref{sec:5_1}. The second model is a real-world reservoir model, and we apply the multiscale approaches to this model to optimize the control strategy. 

\subsubsection{Model 1: 5-spot model}
\label{sec:5_2_1}

\begin{figure}[htbp]
\centering
  \includegraphics[width=0.45\textwidth]{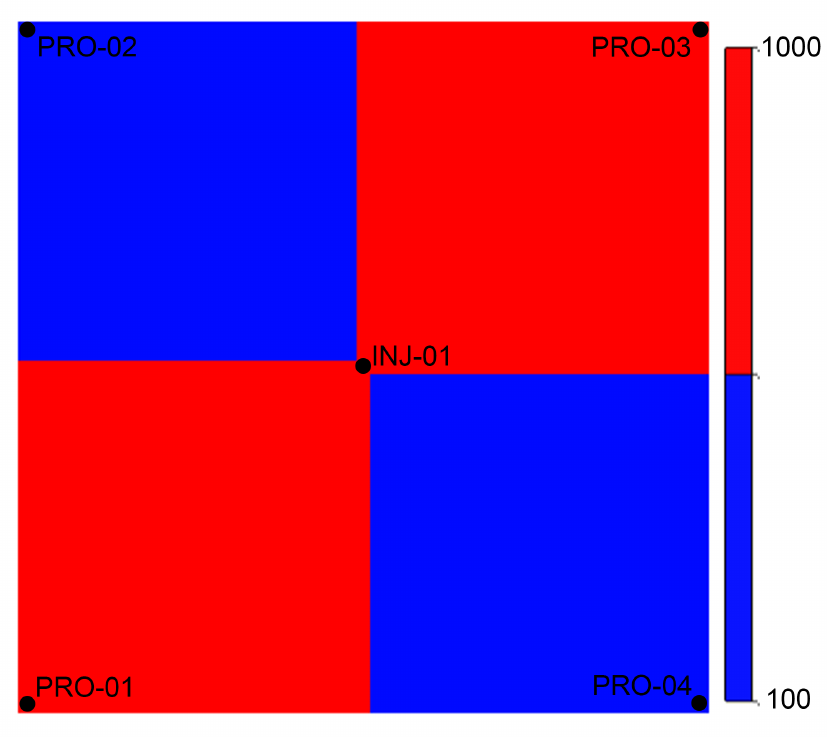}
\caption{The permeanbility field (mD) for model 1.}
\label{fig:perm}
\end{figure}

The first model is a single-layer reservoir containing four producing wells and one injection well in a five-spot well pattern \cite{oliveira_adaptive_2014}. The reservoir model is represented by a 51 $\times$ 51 uniform grid ($\bigtriangleup x=\bigtriangleup y=10 m$; $\bigtriangleup z=5 m$). We consider only oil-water two phase flow. The reservoir permeability field and well placements are shown in Fig. \ref{fig:perm}. We note that there are four different regions of homogeneous permeability. The permeabilities are 1000 mD for the two high-permeability regions, and 100 mD for the two low-permeability regions. The porosity, net-to-gross ratio, and initial water saturation are all 0.2 at all grid blocks. 

The reservoir lifetime is set to 720 days. The injection well INJ-01 (in Fig. \ref{fig:perm}) is not controlled, the liquid rate is fixed at 240 $m^3/d$. The liquid rates of four producing wells are the optimization variables. Bound constraints are considered for the producing wells. The lower bound is set to 0 $m^3/d$ and the upper bound is 80 $m^3/d$ for PRO-01 \& PRO-03 while 0 $m^3/d$ and 40 $m^3/d$ are the lower and upper bounds for PRO-02 \& PRO-04. The initial rates of all producing wells are 20 $m^3/d$.

The objective function we use for this model is the NPV (see equation (\ref{eq:npv})) and the corresponding economic parameters are given in Table \ref{tab:eco_set_1}. 

\begin{table}[htbp]
\centering
\caption{Economic parameters used for model 1.}
\label{tab:eco_set_1}
\begin{tabular}{cc}
\hline\noalign{\smallskip}
Parameter & Value \\
\noalign{\smallskip}\hline\noalign{\smallskip}
Oil revenue & USD 500.0/$m^3$ \\
Water-production cost & USD 250.0/$m^3$ \\
Water-injection cost & USD 80.0/$m^3$ \\
Annual discount rate & 0 \\
\noalign{\smallskip}\hline
\end{tabular}
\end{table}

We use Eclipse 100 \cite{geoquest_eclipse_2014}, a commercial reservoir simulation software from Schlumberger Ltd., to calculate the relevant time-dependent production information for every well for all experiments in this paper.

\subsubsection{Model 2: PUNQ-S3}
\label{sec:5_2_2}

\begin{figure*}[htbp]
\centering 
\subfigure[Tops]{ 
    \label{fig:subfig:punqs3_tops} %% label for first subfigure 
    \includegraphics[width=0.3\textwidth]{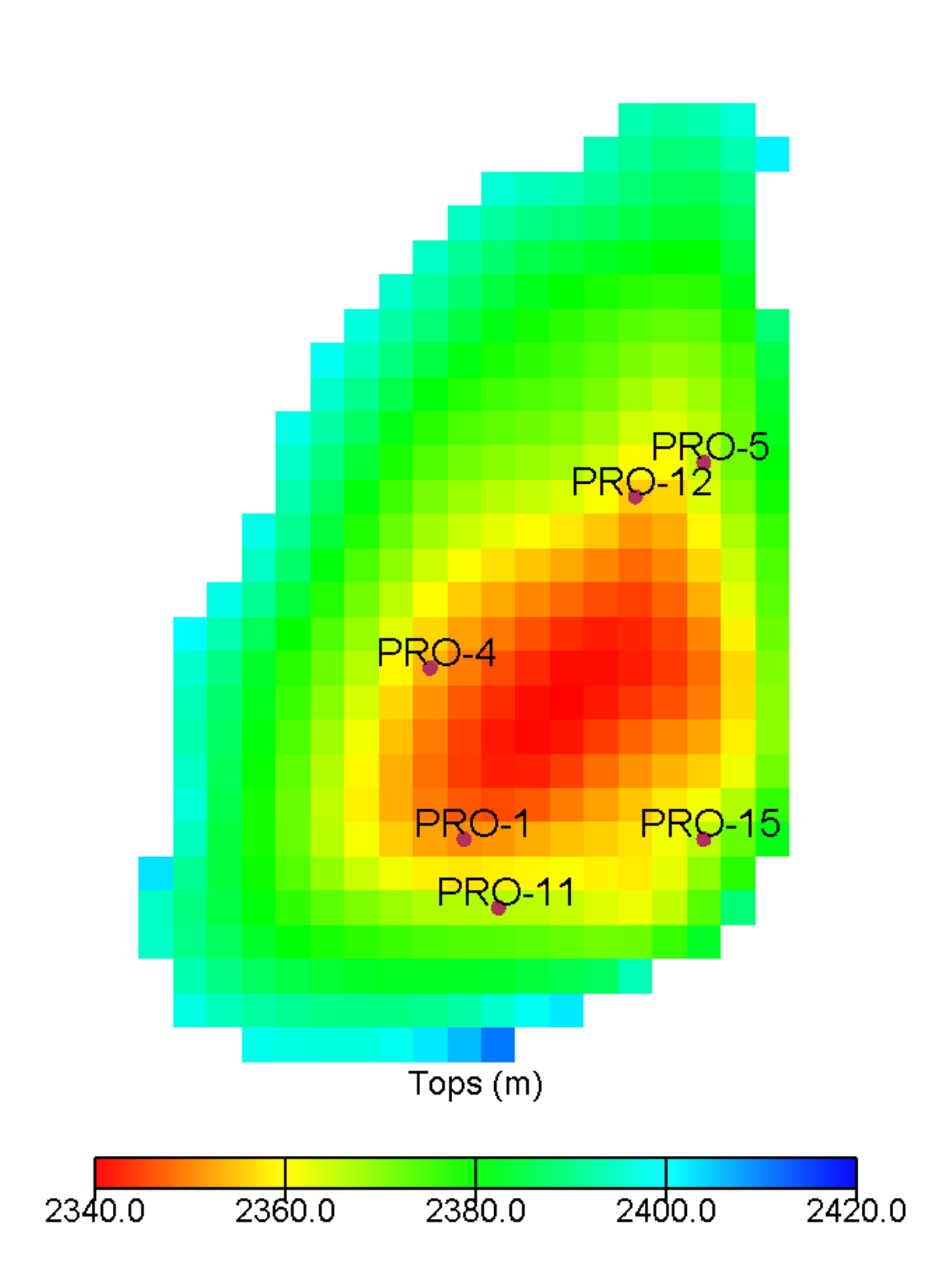}
    }
\subfigure[Permeability]{ 
    \label{fig:subfig:punqs3_perm} %% label for first subfigure 
    \includegraphics[width=0.3\textwidth]{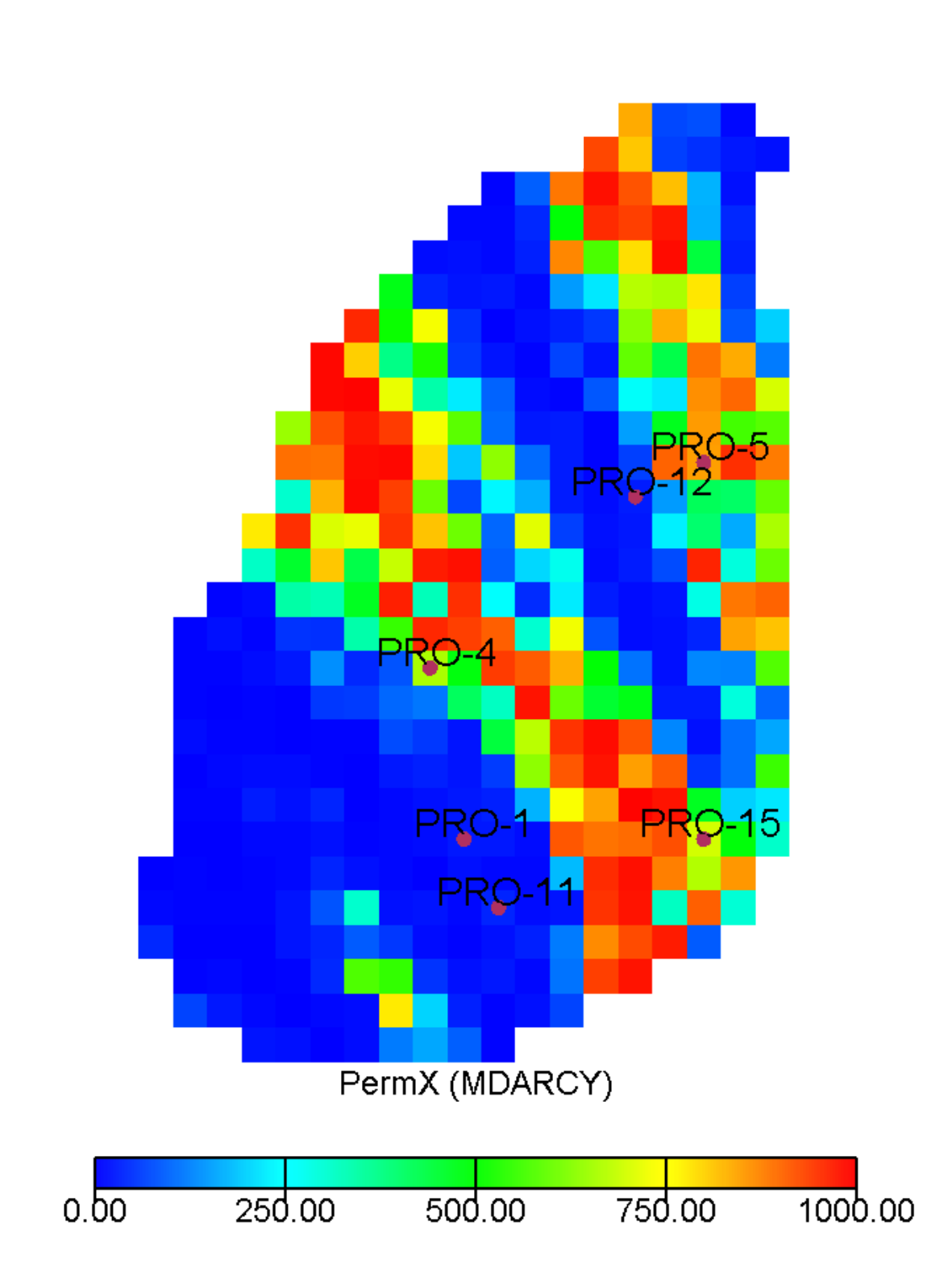}
    }
\subfigure[Oil saturation]{ 
    \label{fig:subfig:punqs3_soil} %% label for first subfigure 
    \includegraphics[width=0.3\textwidth]{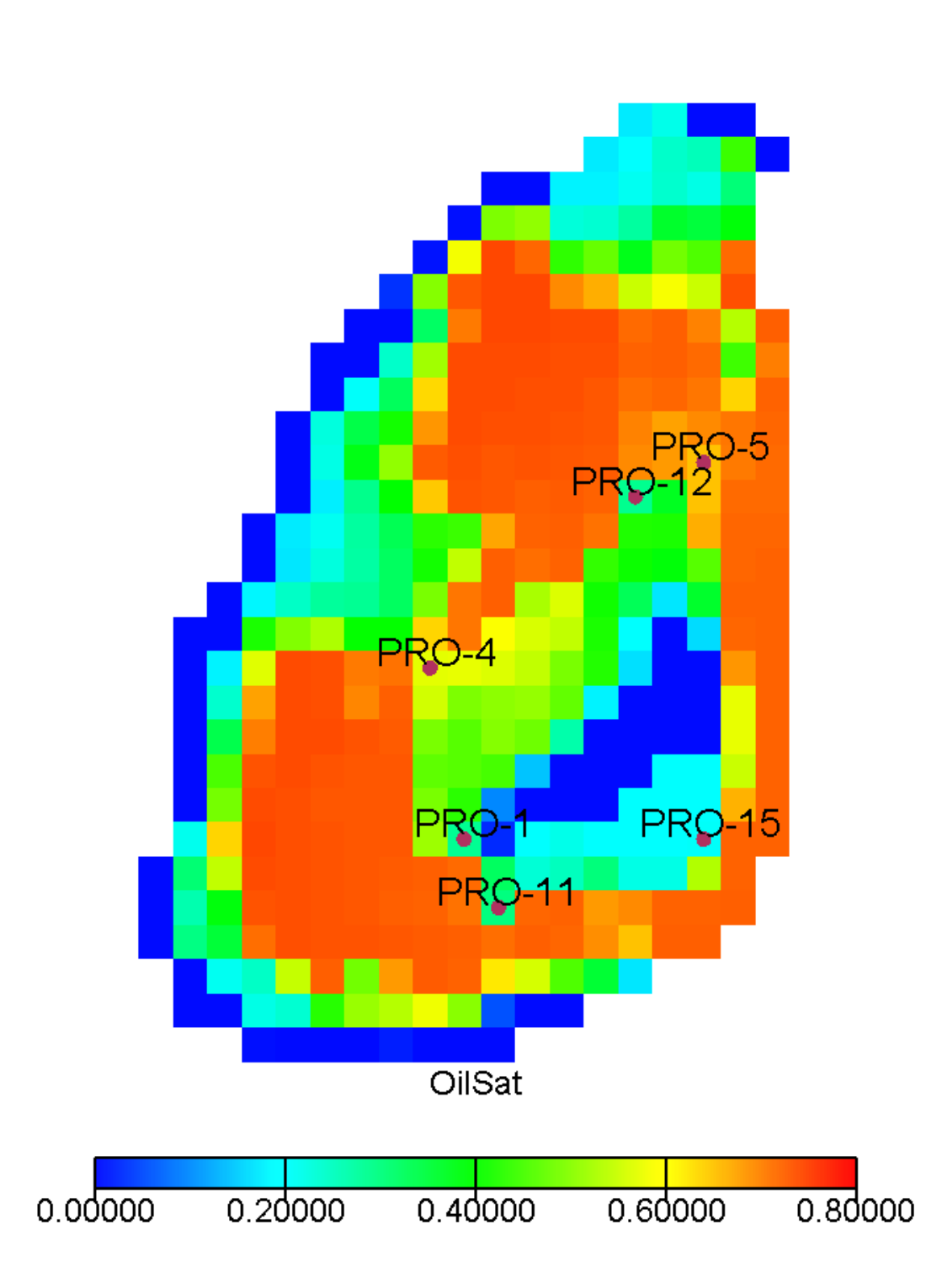}
    }
\caption{Property and wells of PUNQ-S3 field.}
\label{fig:punqs3}
\end{figure*}

The second reservoir model is the PUNQ-S3, which is a small-size reservoir model based on the North Sea reservoir \cite{gao_quantifying_2006}. The model contains a three phase gas-oil-water system with 19 $\times$ 28 $\times$ 5 grid blocks, of which 1761 blocks are active. The field contains 6 production wells but no injection wells are present due to the strong aquifer. Fig. \ref{fig:punqs3} shows the tops (depth of the top phase), permeability and oil saturation present in the model. 

We use a production period of 3840 days (about 10 years), with a minimum control interval of 120 days. The initial liquid rates for all wells are 100 $m^3/d$. The lower bound is set to 0 $m^3/d$ and the upper bound is 200 $m^3/d$ for all wells. BHP bounds are also considered in this example. The lower BHP bounds are set to 12 $\rm{MPa}$ and no upper bound is enforced for any producers. The economic parameters for the NPV calculation are given in Table \ref{tab:eco_set_2}.

\begin{table}[htbp]
\centering
\caption{Economic parameters used for PUNQ-S3.}
\label{tab:eco_set_2}
\begin{tabular}{cc}
\hline\noalign{\smallskip}
Parameter & Value \\
\noalign{\smallskip}\hline\noalign{\smallskip}
Gas revenue & USD 0.5/$m^3$ \\
Oil revenue & USD 500.0/$m^3$ \\
Water-production cost & USD 80.0/$m^3$ \\
Annual discount rate & 0 \\
\noalign{\smallskip}\hline
\end{tabular}
\end{table}

\section{Results and discussion}
\label{sec:6}

\subsection{Effects of control frequency on well control optimization}
\label{sec:6_1}

To show the effect of the control frequency on the NPV for the first model, as described in Section \ref{sec:5_2_1}, we turn off the multiscale approach and optimize using four different, fixed control frequencies. Four control frequencies are considered and these constitute four variations of the optimization problem:

\begin{itemize}
\item[\textbullet] Case 1A, each well is produced under a liquid rate throughout its lifetime. This gives 4 optimization variables in total.
\item[\textbullet] Case 1B, the liquid rate for each well updated every 360 days (2 control periods). This gives 8 optimization variables in total. 
\item[\textbullet] Case 1C, the liquid rate for each well updated every 90 days (8 control periods). This gives 32 optimization variables in total. 
\item[\textbullet] Case 1D, the liquid rate for each well updated every 22.5 days (32 control periods). This gives 128 optimization variables in total. 
\end{itemize}

Three optimization algorithms, GPS, PSO, and CMA-ES, are applied to each case to find the optimal controls and the corresponding NPV.

\subsubsection{Optimal NPV under different control frequencies}
\label{sec:6_1_1}

\begin{figure}[htbp]
\centering
  \includegraphics[width=0.45\textwidth]{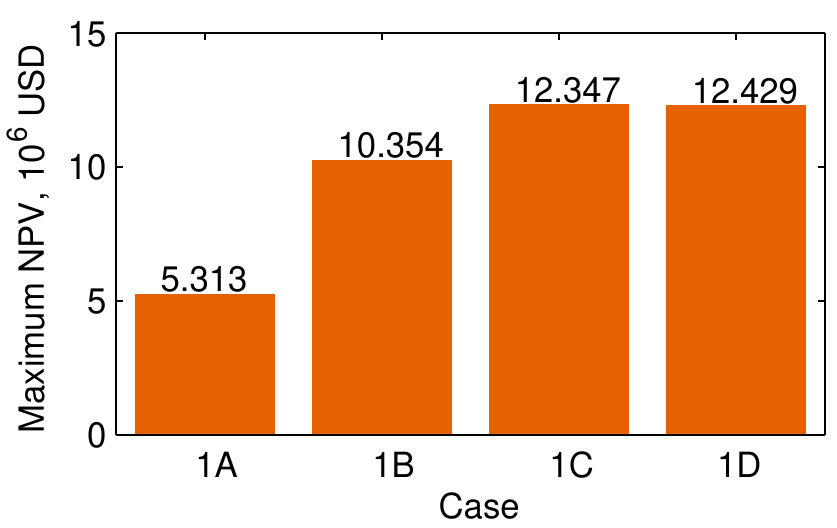}
\caption{Optimum NPV for cases with different control frequencies.}
\label{fig:step_npv}
\end{figure}

Fig. \ref{fig:step_npv} compares the optimal NPV under different control frequencies for this model. The results shown are the best values found using all the optimization approaches. Well control with a reasonable frequency is necessary --- we obtain a significantly higher NPV than what is possible when using a fixed rate over the life cycle (Case 1A). It is clear that with the increase of the number of control adjustments, the optimal NPV grows more and more slowly. The increase in maximum NPV found is very slight (0.28 \%) when the number of control steps for each well increases from 8 (Case 1C) to 32 (Case 1D). There is no need to adjust well rates too frequently. We will not  see a considerable revenue increase and the increase in the number of control adjustments will increase operation costs.  Also the problem with a large  number of control adjustments is harder to optimize and have a higher risk of falling into a local optima (see Section \ref{sec:6_2_1}). 
This justifies the use of multiscale approach to determine an appropriate control frequency.

\subsubsection{Optimal controls under different control frequencies}
\label{sec:6_1_2}

\begin{figure*}[htbp]
\centering 
\subfigure[PRO-01, 1 step]{ 
    \label{fig:subfig:rate1_1} %% label for first subfigure 
    \includegraphics[width=0.22\textwidth]{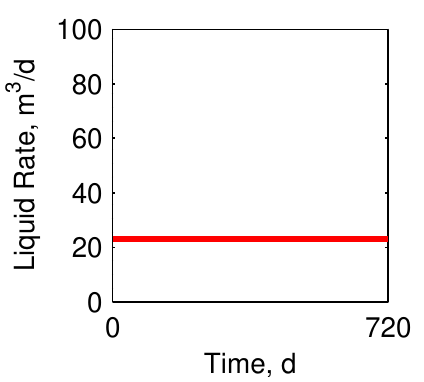}
    }
\subfigure[PRO-01, 2 steps]{ 
    \label{fig:subfig:rate1_2} %% label for first subfigure 
    \includegraphics[width=0.22\textwidth]{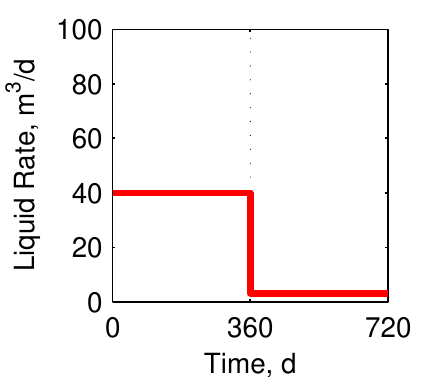}
    }
\subfigure[PRO-01, 8 steps]{ 
    \label{fig:subfig:rate1_3} %% label for first subfigure 
    \includegraphics[width=0.22\textwidth]{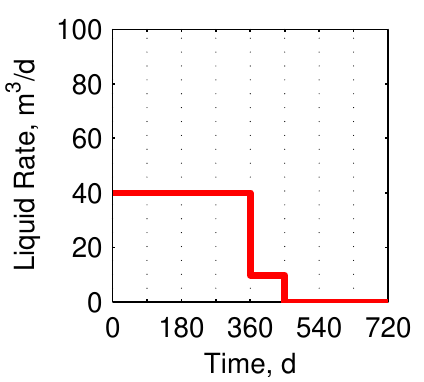}
    }
\subfigure[PRO-01, 32 steps]{ 
    \label{fig:subfig:rate1_4} %% label for first subfigure 
    \includegraphics[width=0.22\textwidth]{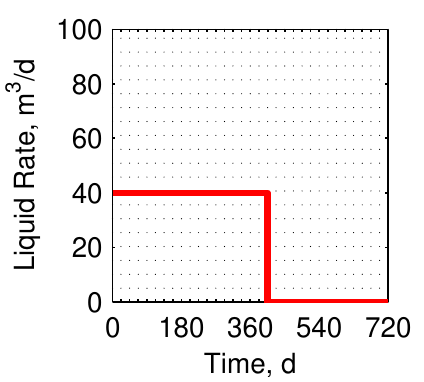}
    }\\
\subfigure[PRO-02, 1 step]{ 
    \label{fig:subfig:rate2_1} %% label for first subfigure 
    \includegraphics[width=0.22\textwidth]{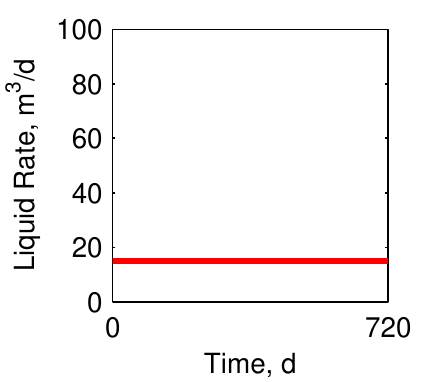}
    }
\subfigure[PRO-02, 2 steps]{ 
    \label{fig:subfig:rate2_2} %% label for first subfigure 
    \includegraphics[width=0.22\textwidth]{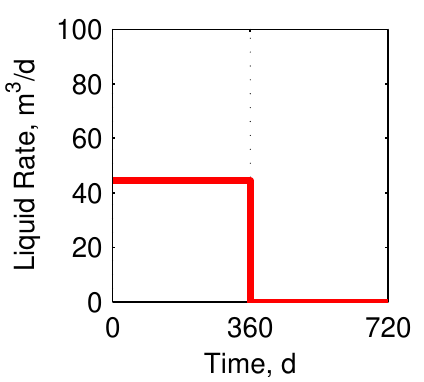}
    }
\subfigure[PRO-02, 8 steps]{ 
    \label{fig:subfig:rate2_3} %% label for first subfigure 
    \includegraphics[width=0.22\textwidth]{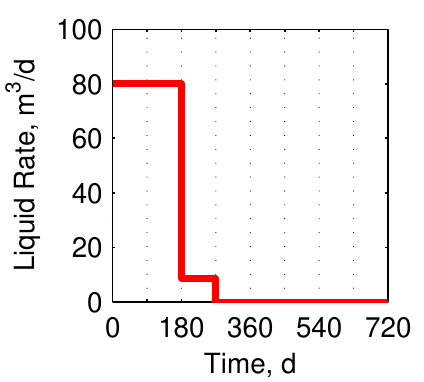}
    }
\subfigure[PRO-02, 32 steps]{ 
    \label{fig:subfig:rate2_4} %% label for first subfigure 
    \includegraphics[width=0.22\textwidth]{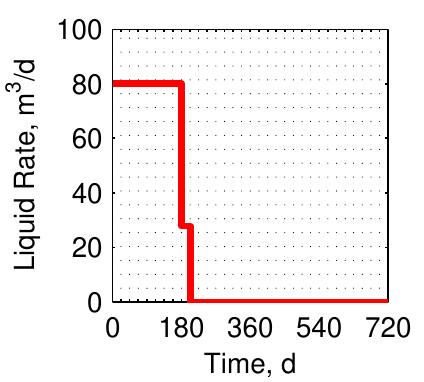}
    }
  \caption{Optimum well control for cases with different control frequencies.} 
  \label{fig:opt_control} %% label for entire figure 
\end{figure*}

Fig. \ref{fig:opt_control} presents optimum controls for wells PRO-01 and PRO-02 under different control frequencies. We omit the results for well PRO-03 and PRO-04 because the reservoir is symmetric. The optimum controls become more like a bang-bang solution for all wells with an increase in the number of control steps. It is worth noting that the optimum controls for Case 1A are significantly different that those for Cases 1B--1D. This reflects the different production strategies for wells using a static rate compared to using dynamic well controls in water flooding reservoirs. The similarity of optimum controls between different control frequencies is important for the success of a multiscale framework. As the multiscale approach uses the optimal controls found in iteration $i$  as the initial guess for iteration $i+1$, a good initial guess could accelerate optimization process and a bad initial guess may mislead the optimizer to wrong search areas and directions (see Section \ref{sec:6_2_3}). Indeed Fig. \ref{fig:opt_control} shows this required similarity as the number of control steps is increased.

\subsection{Performance of GPS, PSO, and CMA-ES for well control optimization}
\label{sec:6_2}

In this section we address the performance of GPS, PSO, and CMA-ES for well control optimization without the use of the multiscale framework. We use the same test cases as in Section \ref{sec:6_1}. We use the maximum number of simulation runs as the stopping criterion, and this value is set to 100 times the number of control variables. As PSO and CMA-ES are stochastic algorithms, 10 trials are performed for these two algorithms to assess the average performance.

\subsubsection{Parameter tuning and the effect of the initial guess on GPS, PSO, and CMA-ES}
\label{sec:6_2_3}

The performance of the optimization algorithms are affected by the choice of their parameter values. In this section, we complete parameter tunings for GPS, PSO, and CMA-ES to improve their performance in solving well control optimization problems.
Here we perform a tuning study for two choices of initial guesses. The {\em good} initial guess is chosen to mimic the initial guess provided by the multiscale algorithm. The {\em bad} initial guess is purposely chosen to be far away from the optimal controls. 
%Additionally, to see the effect of the initial guess we add a {\em bad} quality initial guess as shown in Table \ref{tab:init_guess}.
%To prepare to hybridize GPS, PSO, and CMA-ES with the multiscale framework, we mainly focus on the performance of these algorithms with different initial guesses since the key idea of the multiscale approach is to use the coarse-scale solution as the initial guess for the next finer scale optimization.

We hypothesize that the performance of the local search algorithms are highly affected by the initial guess, while the stochastic global search algorithms are not. We take Case 1B as an example and use the three different initial guesses shown in Table \ref{tab:init_guess}. For each initial guess, 10 trials were performed for PSO and CMA-ES and 1 trial for GPS (since it is a deterministic algorithm). 

\begin{table}[htbp]
\centering
\caption{Three initial guesses for GPS, PSO, and CMA-ES.}
\label{tab:init_guess}
\begin{tabular}{ccc}
\hline
Type & Initial guess & NPV, $\times 10^6$ USD \\
\hline
good   & $[20,20,20,20,20,20,20,20]$ & 5.0009  \\
medium & $[20,20,40,40,40,40,20,20]$ & 2.6484 \\
bad    & $[ 0,40, 0,80, 0,80, 0,40]$ & -4.2826 \\
\hline
\end{tabular}
\end{table}

\begin{figure}[htbp]
\centering 
\subfigure[GPS]{ 
    \label{fig:subfig:i_gps} %% label for first subfigure 
    \includegraphics[width=0.4\textwidth]{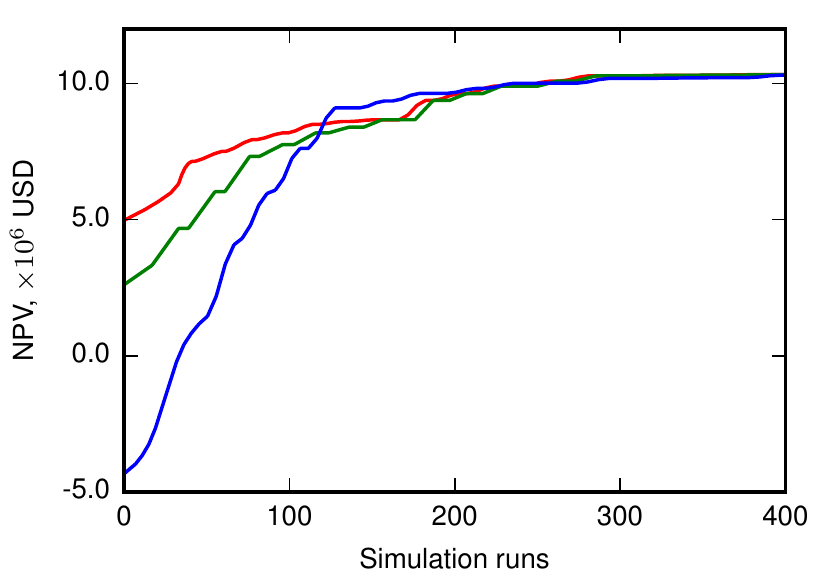}
    }
\subfigure[PSO]{ 
    \label{fig:subfig:i_pso} %% label for first subfigure 
    \includegraphics[width=0.4\textwidth]{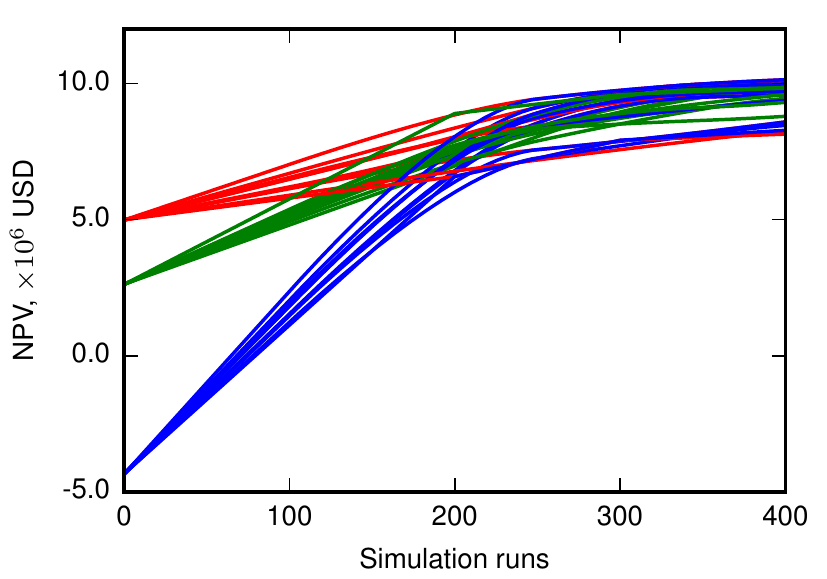}
    }
\subfigure[CMA-ES]{ 
    \label{fig:subfig:i_cmaes} %% label for first subfigure 
    \includegraphics[width=0.4\textwidth]{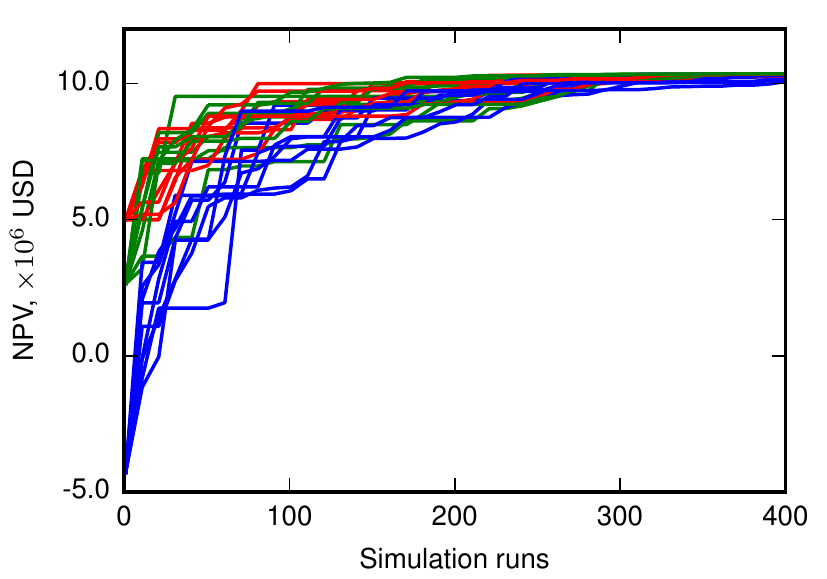}
    }
  \caption{Plots of NPV versus the number of simulation runs for GPS, PSO, and CMA-ES with different initial guesses with 8 control steps. Each line represents one trial. The color red, green, blue denotes the good, medium, bad initial guess, respectively.} 
  \label{fig:case2_init} %% label for entire figure 
\end{figure}

Fig. \ref{fig:case2_init} shows the plots of NPV versus the number of simulations for GPS, PSO, and CMA-ES. Each line represents one trial. 
The early convergence of GPS and PSO is affected by the choice of the initial guess. CMA-ES recovers quite quickly from the bad initial guess, even converges more quickly than from a good initial guess in some cases.
%It is clear that, for GPS and PSO, with the different initial guesses, the performance of these two algorithms differed significantly. For CMA-ES, the difference is not so obvious. CMA-ES with a bad initial guess converges faster than the one with better guess in some cases. 
With a large number of simulation runs, the effect of initial guess for all three algorithms is quite small. This suggests that if we want to make efficient use of the multiscale approach, the number of simulations at each scale should be limited.

Again using a good and a bad initial guess, we analyzed the effect of the other parameters for PSO and CMA-ES, to find out the primary parameters which affect the performance of the algorithms. For PSO, the algorithm parameters include the population size $\lambda$, and the parameters $\omega$, $c_1$, and $c_2$. Three levels are chosen for each parameter. For CMA-ES, the algorithm parameters include the population size $\lambda$, the parent number $\mu$ (number of candidate solutions used to update the distribution parameters), the recombination weights $\omega$, and the parameter $\sigma$, which determines the initial coordinate-wise standard deviations for the search. 

\begin{figure*}[htbp]
\centering 
\subfigure[PSO, $\lambda$]{ 
    \label{fig:subfig:pt_pso_lambda} %% label for first subfigure 
    \includegraphics[width=0.4\textwidth]{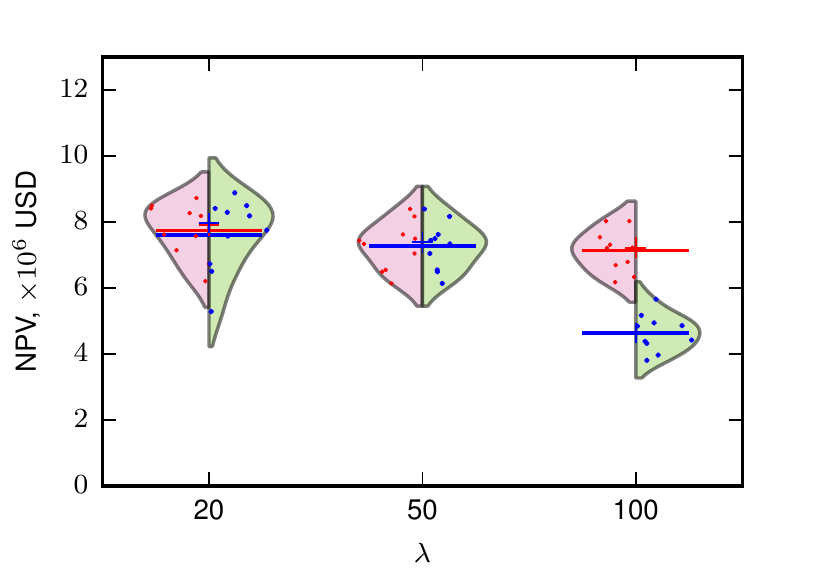}
    }
\subfigure[PSO, $\omega$]{ 
    \label{fig:subfig:pt_pso_omega} %% label for first subfigure 
    \includegraphics[width=0.4\textwidth]{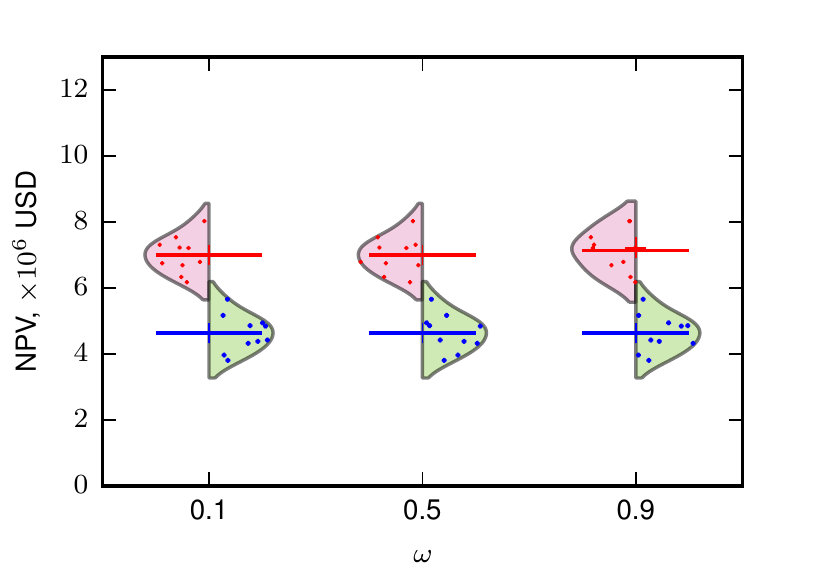}
    }
\subfigure[PSO, $c_1$]{ 
    \label{fig:subfig:pt_pso_c1} %% label for first subfigure 
    \includegraphics[width=0.4\textwidth]{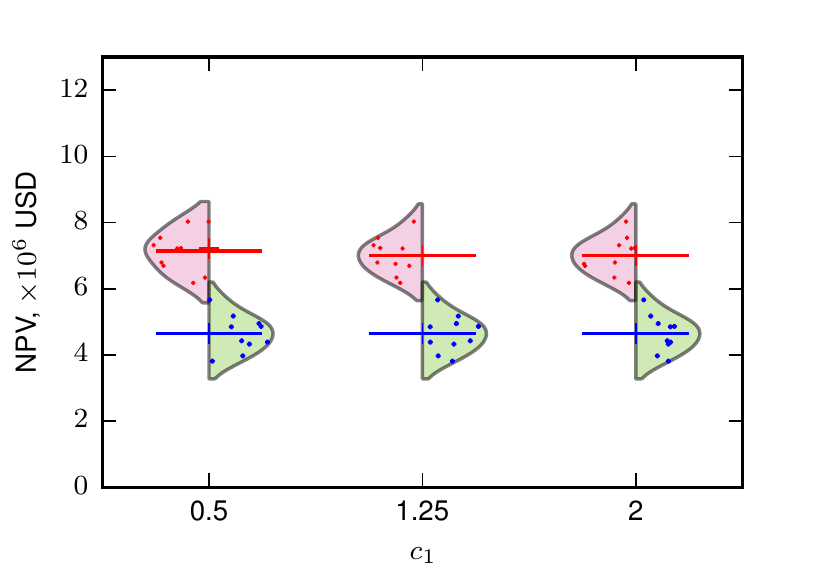}
    }
\subfigure[PSO, $c_2$]{ 
    \label{fig:subfig:pt_pso_c2} %% label for first subfigure 
    \includegraphics[width=0.4\textwidth]{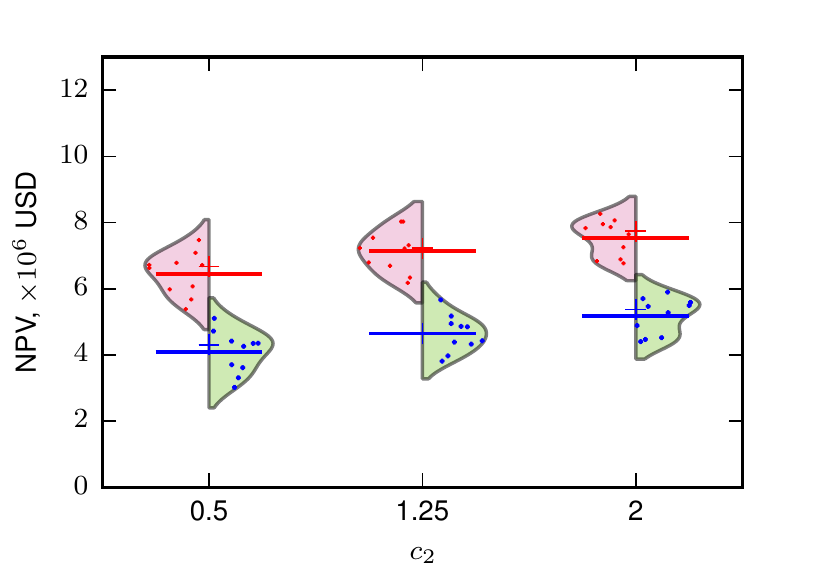}
    }
\caption{Beanplots of the NPV for various parameter settings of PSO. The left side of each beanplot gives the results obtained with a good initial guess, and the right side gives the results obtained with a bad initial guess. The individual dots show the NPV obtained by each trial. The background pink and green colors show the distribution of results. The short horizontal line and the marker `+' denote the mean and median of all 10 trials, respectively.} 
\label{fig:pso_paratuning} %% label for entire figure 
\end{figure*}

\begin{figure*}[htbp]
\centering 
\subfigure[CMA-ES, $\sigma$]{ 
    \label{fig:subfig:pt_cmaes_sigma} %% label for first subfigure 
    \includegraphics[width=0.4\textwidth]{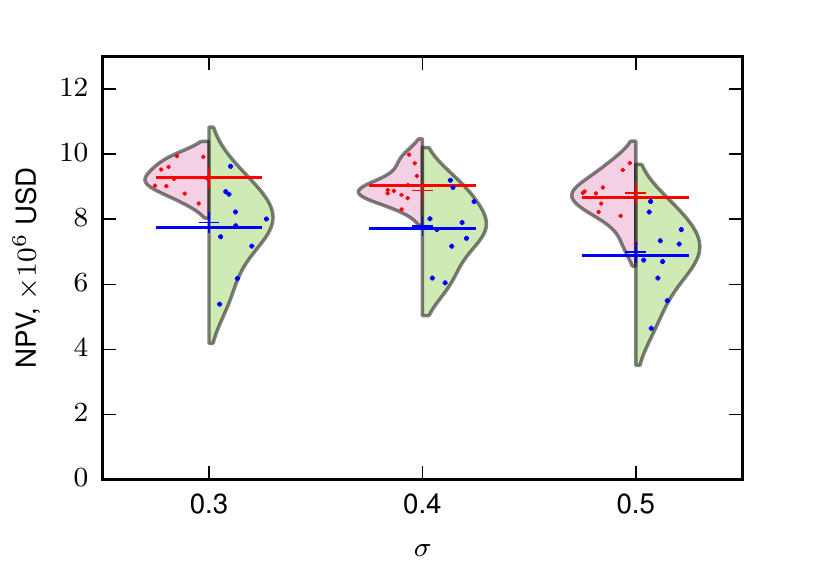}
    }
\subfigure[CMA-ES, $\lambda$]{ 
    \label{fig:subfig:pt_cmaes_lambda} %% label for first subfigure 
    \includegraphics[width=0.4\textwidth]{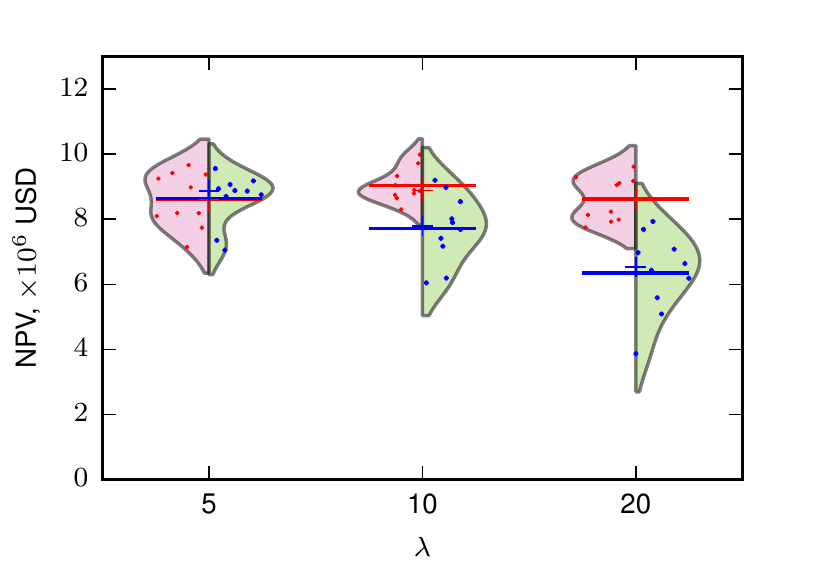}
    }
\subfigure[CMA-ES, $\mu$]{ 
    \label{fig:subfig:pt_cmaes_mu} %% label for first subfigure 
    \includegraphics[width=0.4\textwidth]{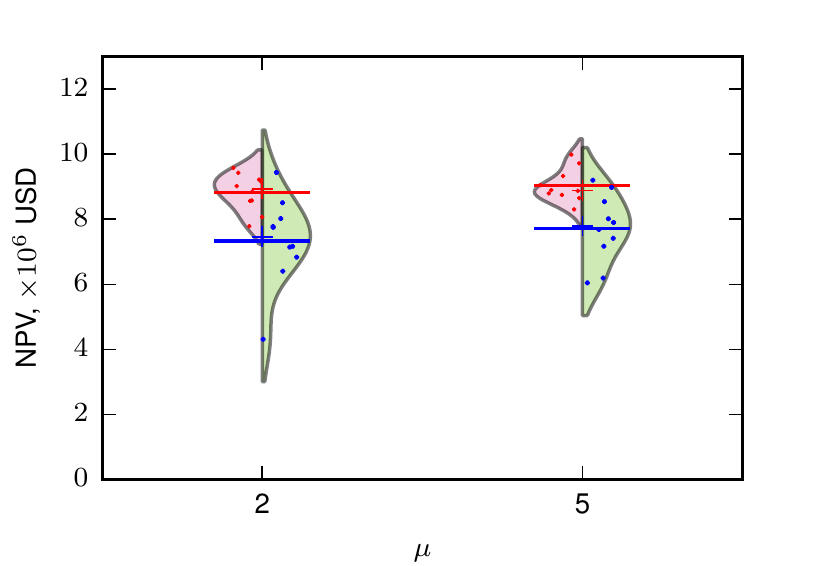}
    }
\subfigure[CMA-ES, $\omega$]{ 
    \label{fig:subfig:pt_cmaes_omega} %% label for first subfigure 
    \includegraphics[width=0.4\textwidth]{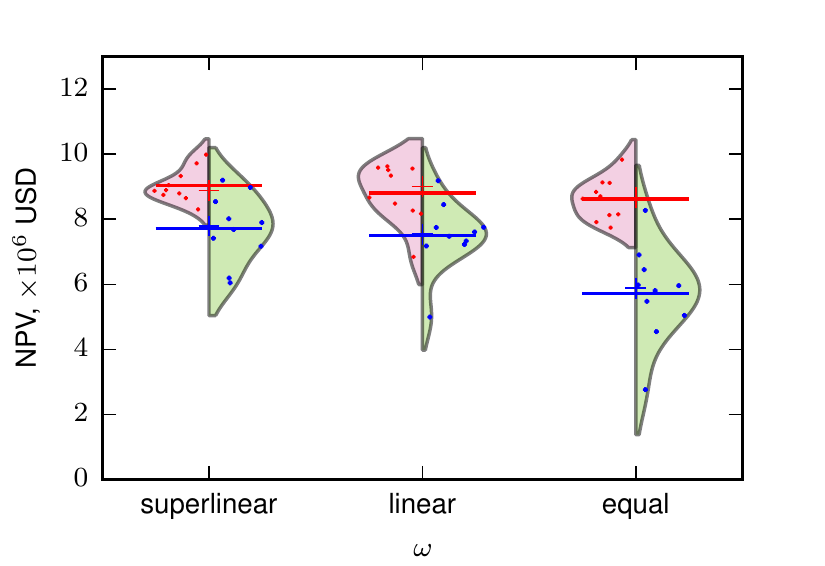}
    }
\caption{Beanplots of the NPV for various parameter settings of CMA-ES. The left side of each beanplot gives the results obtained with a good initial guess, and the right side gives the results obtained with a bad initial guess. The individual dots show the NPV obtained by each trial. The background pink and green colors show the distribution of results. The short horizontal line and the marker `+' denote the mean and median of all 10 trials, respectively.} 
\label{fig:cmaes_paratuning} %% label for entire figure 
\end{figure*}

For each parameter choice, 10 trials are performed for Case 1B with a good and a bad initial guess to start the optimization. We use the best NPV obtained after 20\% of maximum simulation runs as the evaluation criterion of the algorithm's performance. We mainly focus on the performance of the algorithms at an early stage because the multiscale framework requires a good early stage performance for the hybrid optimization algorithm. Also our test results showed that the performance of the algorithms are less sensitive to the parameter values at a later stage. The performance for each parameter choice is shown in the beanplots in Fig. \ref{fig:pso_paratuning} and Fig. \ref{fig:cmaes_paratuning}. A beanplot \cite{kampstra_beanplot:_2008} promotes visual comparison of univariate data between groups. In a beanplot, the individual observations are shown as small points or small lines in a one-dimensional scatter plot. In addition, the estimated density of the distributions is visible and the mean (bold line) and median (marker `+') are shown. 

From Fig. \ref{fig:pso_paratuning} we can see that for PSO the population size plays the most important role in the algorithm's ability to utilize the good initial guess. When the population size equals 20 or 50, PSO with a bad initial guess obtains a similar NPV as PSO with a good initial guess. This is because with the same number of simulation runs, PSO with a small population size can evolve more generations, and this decreases the affect of the initial guess. 
%The bigger population size, the smaller variation for the algorithm. 
The bigger population size, the smaller the variability in the NPV results with a similar mean value. For these reasons, we choose $\lambda=100$ for all subsequent PSO experiment.
Parameter $c_2$ is one of the weighting parameters, the bigger $c_2$, the greater the tendency for the particles to fly towards the best location found so far. We suggest a bigger $c_2$ when combining with the multiscale approach. The parameters $\omega$ and $c_1$ have no obvious affect in this case. Generally for all parameter values PSO responds favorably to the better initial guess.

From Fig. \ref{fig:cmaes_paratuning} we can see that for CMA-ES the good initial guess gives a higher average NPV with smaller variability.
%all four parameters affect the role of the initial guess. 
For this problem, the best configuration is $\sigma=0.3$, $\lambda=10$, $\mu=2$, and $\omega=\rm{superlinear}$; this is also the default configuration of CMA-ES. In fact, according to the work of \cite{hansen_benchmarking_2010}, since finding good parameters is considered as part of the algorithm design, CMA-ES does not require significant parameter tuning for its application. The choice of parameters is not left to the user (arguably with the exception of population size $\sigma$).

\subsubsection{Performance with different control frequencies}
\label{sec:6_2_1}

\begin{figure*}[htbp]
  \centering 
  \subfigure[Case 1A]{ 
    \label{fig:subfig:c_s_1a} %% label for first subfigure 
    \includegraphics[width=0.4\textwidth]{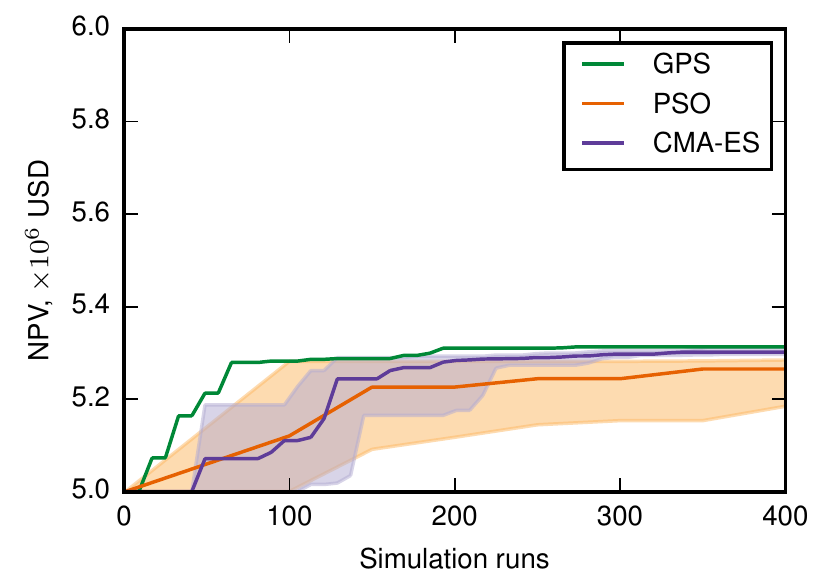}} 
  %\hspace{1in} 
  \subfigure[Case 1B]{ 
    \label{fig:subfig:c_s_1b} %% label for second subfigure 
    \includegraphics[width=0.4\textwidth]{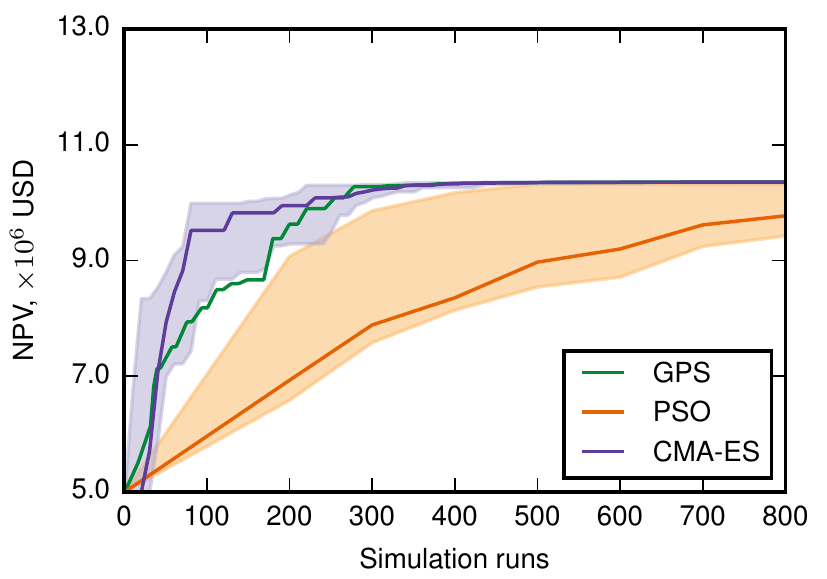}} 
  \subfigure[Case 1C]{ 
    \label{fig:subfig:c_s_1c} %% label for first subfigure 
    \includegraphics[width=0.4\textwidth]{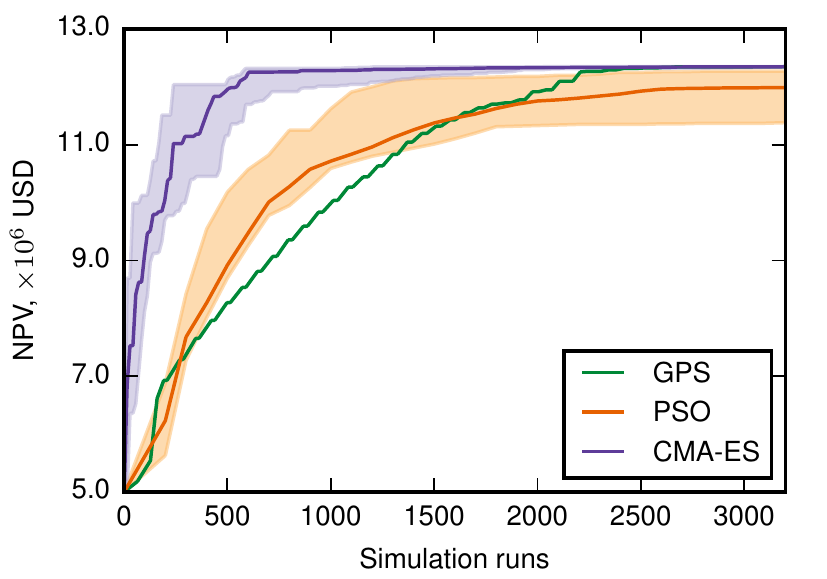}} 
  %\hspace{1in} 
  \subfigure[Case 1D]{ 
    \label{fig:subfig:c_s_1d} %% label for second subfigure 
    \includegraphics[width=0.4\textwidth]{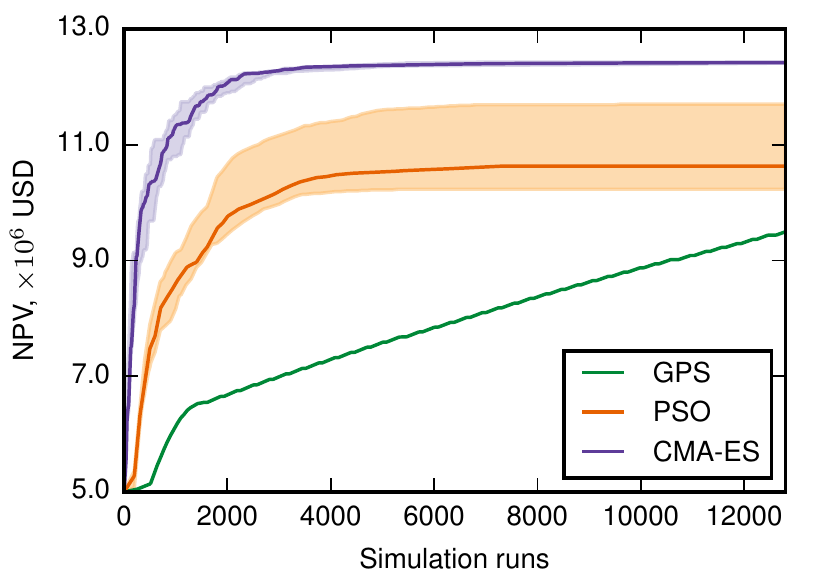}}   
  \caption{Optimization performance for well control problems using GPS, PSO and CMA-ES. The solid lines are median NPV over all 10 runs of PSO and CMA-ES without the multiscale framework. And the areas between maximum and minimum NPV are filled with the corresponding color. Note that the $x$-axis scale is different for each case.} 
  \label{fig:c_s} %% label for entire figure 
\end{figure*}

Table \ref{tab:e1_result} and Fig. \ref{fig:c_s} show the performance of GPS, PSO, and CMA-ES for well control optimization problems with different control frequencies. In Table \ref{tab:e1_result}, the maximum, minimum, mean, median, and standard deviation of the NPV for each case are given. From the table we can see that, for Case 1A, which has only 4 variables, GPS obtains the highest NPV after 400 simulation runs. Similar results are found in Case 1B. In Case 1C, the maximum NPV of CMA-ES exceeds the GPS, but the mean and median NPV for CMA-ES are lower than these of GPS. In these three cases, although the final NPV of GPS is larger than the final NPV for CMA-ES and PSO, the difference of the mean/median NPV for the three algorithms is quite small (less than 2\%). In Case 1D, which has 128 variables, the NPV obtained by GPS is obviously lower than CMA-ES. Generally, CMA-ES showed excellent performance in most cases. GPS performs best when the problem dimension is very small.

\begin{table*}[htbp]
\centering
\caption{Results for Cases 1A-1D for Model 1 using GPS, PSO, and CMA-ES. Values shown are NPV ($\times 10^6$ USD). }
\label{tab:e1_result}
\begin{tabular}{llrrrrrr}
\hline
Case & Algorithm & Trials & Max & Min & Mean & Median & Std. \\
\hline
1A & GPS    & 1  &5.3132 &5.3132 &5.3132 &5.3132 &-     \\
   & PSO    & 10 &5.2850 &5.1850 &5.2603 &5.2720 &0.0310 \\
   & CMA-ES & 10 &5.3121 &5.2969 &5.3034 &5.3045 &0.0048 \\
\hline
1B & GPS    & 1  &10.3539 &10.3539 &10.3539 &10.3539 &- \\
   & PSO    & 10 &10.3200 & 9.4220 &10.0840 &10.1700 &0.2819 \\
   & CMA-ES & 10 &10.3536 &10.3511 &10.3527 &10.3528 &0.0008 \\
\hline
1C & GPS    & 1  &12.3470 &12.3470 &12.3470 &12.3470 &- \\
   & PSO    & 10 &12.2700 &11.3800 &11.9660 &12.1050 &0.2966 \\
   & CMA-ES & 10 &12.3474 &12.3447 &12.3466 &12.3467 &0.0007 \\
\hline
1D & GPS    & 1  & 9.5083 & 9.5083 & 9.5083 & 9.5083 &- \\
   & PSO    & 10 &11.7000 &10.2300 &11.1290 &11.1700 &0.4910 \\
   & CMA-ES & 10 &12.4285 &12.3466 &12.4054 &12.4178 &0.0315 \\
\hline
\end{tabular}
\end{table*}

Fig. \ref{fig:c_s} shows the plots of NPV versus the number of simulation runs for the four cases. In this figure, we use a solid line to show the median NPV of each algorithm, and use the same color as the line to fill the area between the maximum and minimum NPV for each algorithm. These plots clearly show the performance of GPS, PSO, and CMA-ES using different computational budgets. GPS obtains the highest NPV for Case 1A--1C at the end of optimization. But the budget (number of simulation runs) required for GPS grows rapidly as the dimension of the problem increases. GPS converged with no more than 50\% of total budget for Cases 1A and 1B, and about 80\% of the total budget for Case 1C. For Case 1D, GPS did not converge after $100D$ simulation runs. CMA-ES obtains almost as high a NPV as GPS for Case 1A--1C, and it obtains highest NPV for Case 1D. Furthermore, CMA-ES showed an excellent performance when the budget is limited. PSO outperforms GPS for a low budget and a large problem dimension, but it still can not beat CMA-ES in these cases.

Since PSO and CMA-ES are stochastic algorithms, the performance is different for each trial. In Table \ref{tab:e1_result} we can see the standard deviation for PSO is larger than the standard deviation for CMA-ES. In Figure \ref{fig:c_s} we can see that the best NPV obtained has a higher variation for low computational budgets than for high budgets. For PSO, the variability did not decrease in Case 1C and 1D as the algorithm converged. 
%the range of NPV between the trials increases first and then decreases for CMA-ES. But for PSO, this range did not decreases in Case 1C and 1D. 

\begin{figure*}[htbp]
  \centering 
  \subfigure{ 
  \label{fig:subfig:indi_lengend}
    \includegraphics[width=0.9\textwidth]{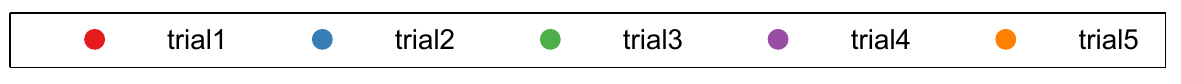}} 
  %\hspace{1in} 
  \setcounter{subfigure}{0}
  \subfigure[CMA-ES (1)]{ 
  \label{fig:subfig:indi_cma_1}
    \includegraphics[width=0.22\textwidth]{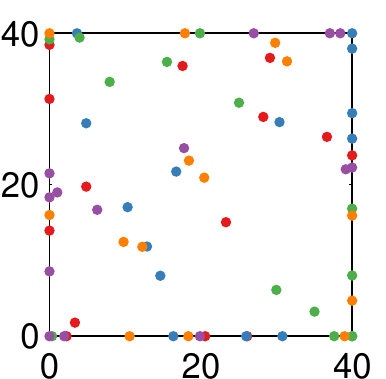}} 
   \subfigure[CMA-ES (36)]{ 
   \label{fig:subfig:indi_cma_2}
    \includegraphics[width=0.22\textwidth]{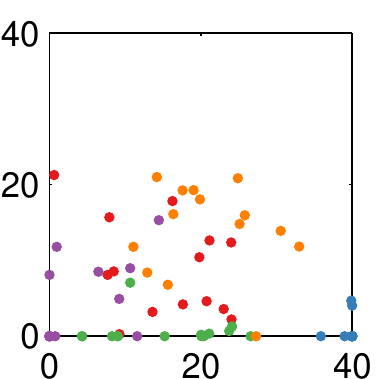}} 
    \subfigure[CMA-ES (72)]{ 
    \label{fig:subfig:indi_cma_3}
    \includegraphics[width=0.22\textwidth]{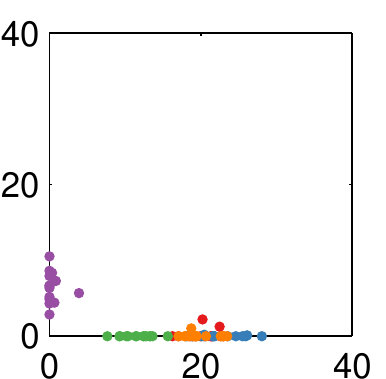}} 
    \subfigure[CMA-ES (215)]{ 
    \label{fig:subfig:indi_cma_4}
    \includegraphics[width=0.22\textwidth]{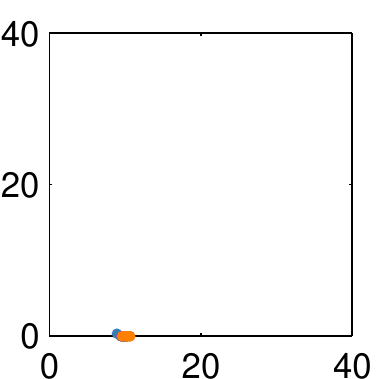}} 
    
      \subfigure[PSO (1)]{ 
  \label{fig:subfig:indi_pso_1}
    \includegraphics[width=0.22\textwidth]{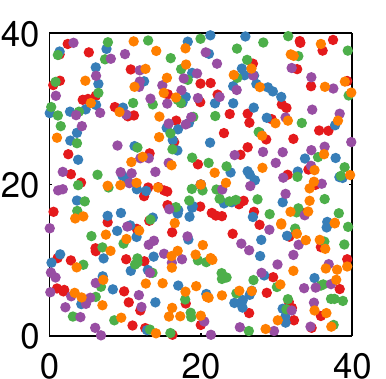}} 
   \subfigure[PSO (5)]{ 
   \label{fig:subfig:indi_pso_2}
    \includegraphics[width=0.22\textwidth]{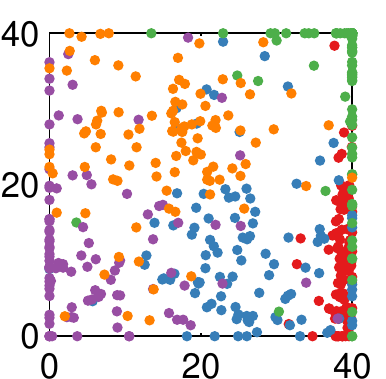}} 
    \subfigure[PSO (10)]{ 
    \label{fig:subfig:indi_pso_3}
    \includegraphics[width=0.22\textwidth]{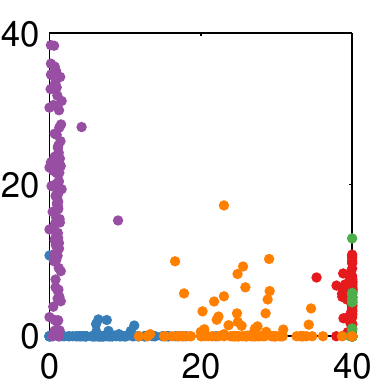}} 
    \subfigure[PSO (30)]{ 
    \label{fig:subfig:indi_pso_4}
    \includegraphics[width=0.22\textwidth]{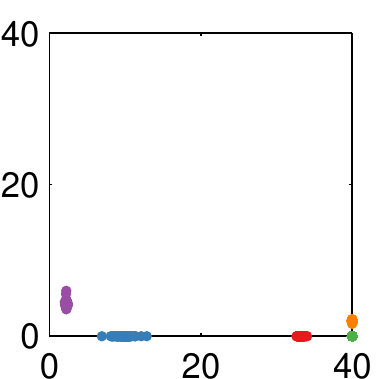}} 
  \caption{Scatter diagrams for CMA-ES and PSO for Case 1C. The bracketed number in the caption of each sub-figure is the 
  iteration number. Points represent the candidate solutions at this specific iteration.} 
  \label{fig:indi_cmapso} %% label for entire figure 
\end{figure*}

To investigate further, we choose 2 of the 32 variables and 5 of the 10 trials for Case 1C and then compare the population distribution of CMA-ES and PSO at different iterations. The resulting scatter diagrams are shown in Fig. \ref{fig:indi_cmapso}. In Case 1C, the population size is 100 for PSO, and 14 for CMA-ES. Hence after same number of simulation runs, CMA-ES and PSO are at different iteration numbers. After 500 simulation runs, CMA-ES is at the 36th iteration, while PSO is at the 5th iteration. After about 3000 simulation runs (Fig. \ref{fig:subfig:indi_cma_4} and \ref{fig:subfig:indi_pso_4}), we can see that PSO has converged to different locations for each trial. Compared with CMA-ES, PSO is more easily falls into local optima in our test cases, in spite of the larger population size and the ability search the entire space.

\subsubsection{Performance in parallel environments}
\label{sec:6_2_2}

GPS, PSO, and CMA-ES parallelize naturally. Here we investigate the performance of these three algorithms in parallel environments. Table \ref{tab:e1_num_s} gives the population sizes of GPS, PSO, and CMA-ES in Cases 1A--1D. The population size for GPS and CMA-ES are decided automatically by the algorithms based on the problem dimension. For a $D$-dimensional problem, the population size equals $2D$ for GPS, and $4+\lfloor 3\ln{(D)}\rfloor$ for CMA-ES. The population size for PSO is usually decided by the user and we use 100 for all cases (more discussion on the population size is given in Section \ref{sec:6_2_3}). 

In a parallel environment, we can evaluate a number of individuals, up to the number of processors, simultaneously. Note that we are not able to evaluate the individuals from different iterations at the same time. For well control optimization problems, the computation time is mainly spent evaluating the reservoir simulation, a parallel environment can greatly reduce the time of optimization. 

\begin{table}[htbp]
\centering 
\caption{Population size of GPS, PSO, and CMA-ES for Case 1A-1D. $D$ is the number of variables in the problem.}
\label{tab:e1_num_s}
\begin{tabular}{ccccc}
\hline
Case & $D$ & GPS & PSO & CMA-ES \\
\hline
1A & 4   & 8   & 100 & 8 \\
1B & 8   & 16  & 100 & 10\\
1C & 32  & 64  & 100 & 14\\
1D & 128 & 256 & 100 & 18\\
\hline
\end{tabular}
\end{table}

\begin{figure*}[htbp]
\centering 
\subfigure[Case 1A]{ 
    \label{fig:subfig:1a_p} %% label for first subfigure 
    \includegraphics[width=0.4\textwidth]{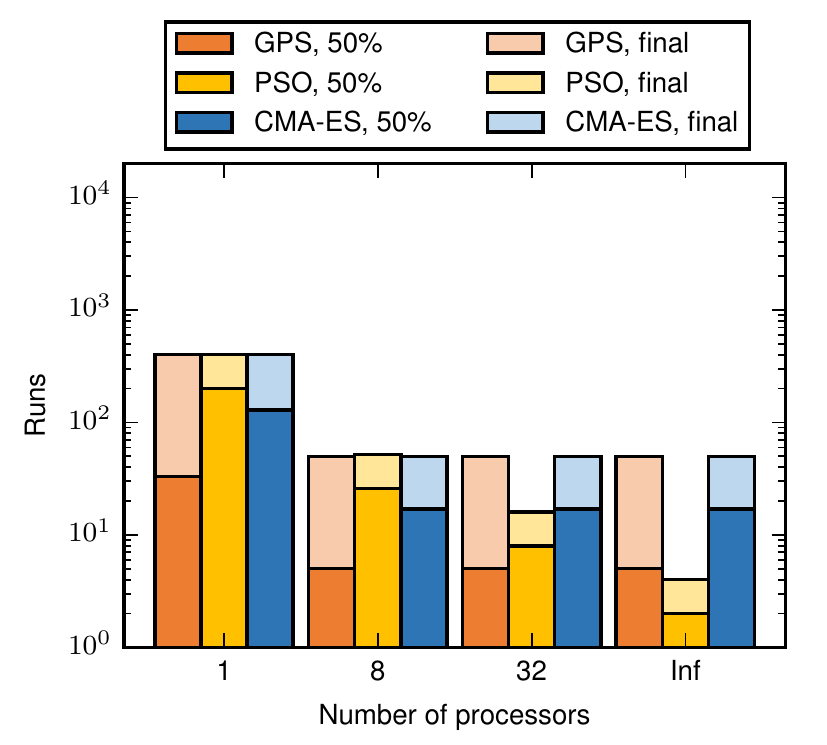}
    }
\subfigure[Case 1B]{ 
    \label{fig:subfig:1b_p} %% label for first subfigure 
    \includegraphics[width=0.4\textwidth]{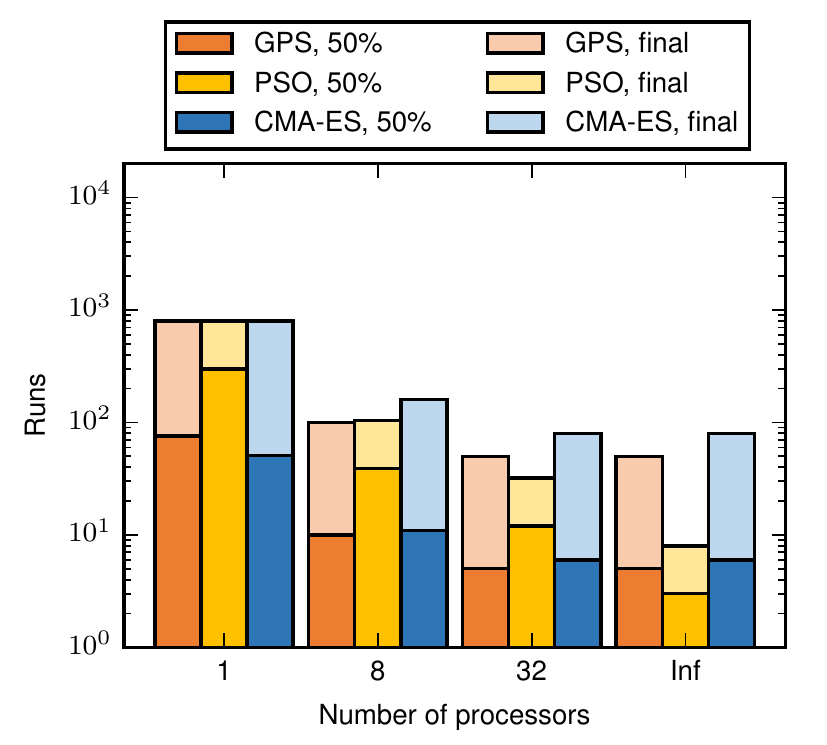}
    }
\\
\subfigure[Case 1C]{ 
    \label{fig:subfig:1c_p} %% label for first subfigure 
    \includegraphics[width=0.4\textwidth]{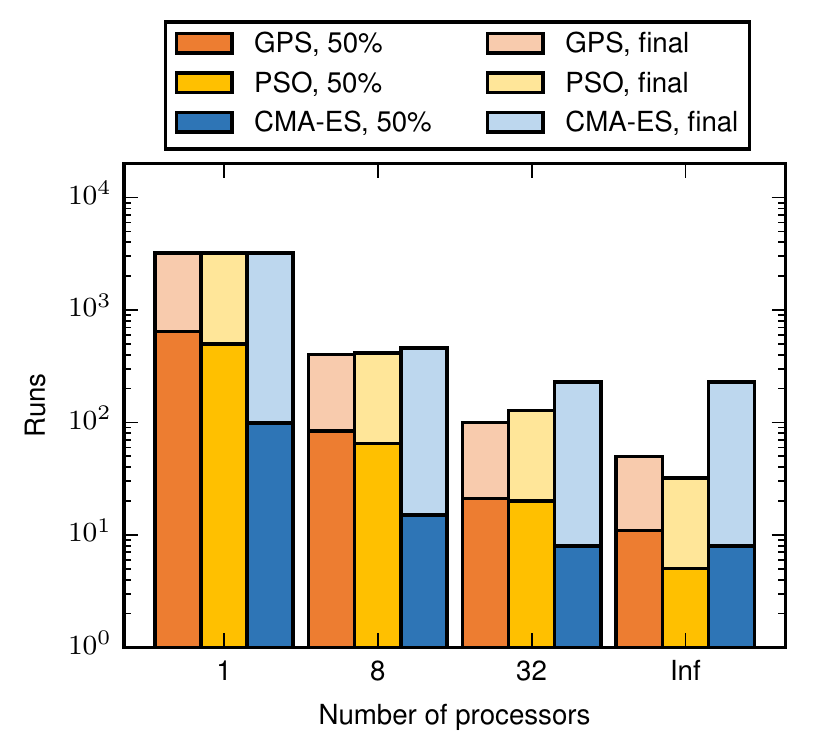}
    }
\subfigure[Case 1D]{ 
    \label{fig:subfig:1d_p} %% label for first subfigure 
    \includegraphics[width=0.4\textwidth]{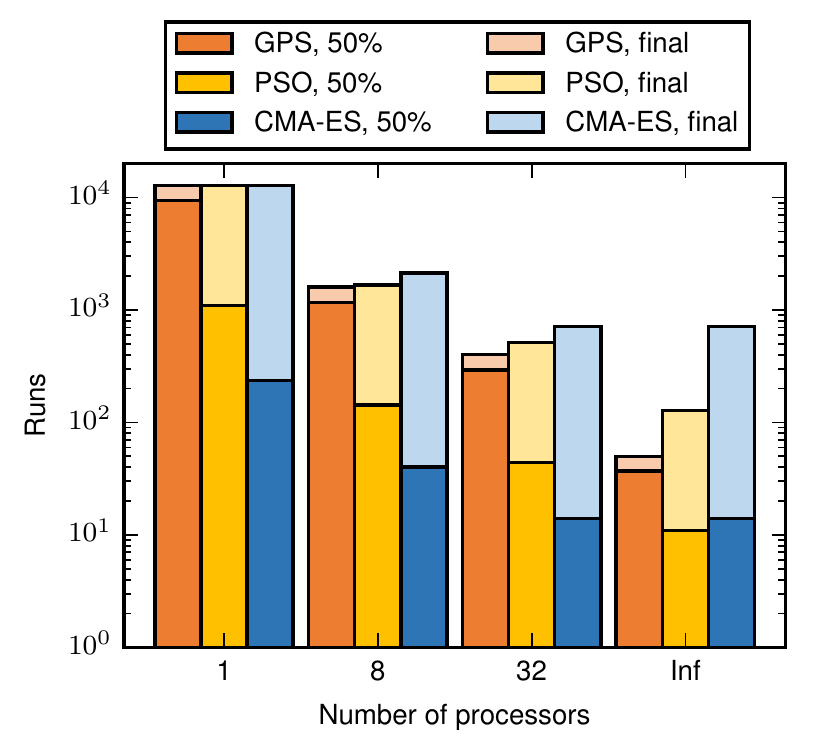}
    } 
  \caption{Number of runs required for the well control optimization problems in parallel environments with different numbers of processors. ``50\%'' in the legend denotes the number of runs required to reach 50\% of the maximum NPV for the algorithms. ``final'' in the legend denotes the total number of runs required to reach the maximum number of simulator evaluation for the algorithms.} 
  \label{fig:case1_p} %% label for entire figure 
\end{figure*}

Assume we have three parallel environments, with 8, 32, and an infinite number of processors, respectively. Fig. \ref{fig:case1_p} compares the parallel performance of GPS, PSO, and CMA-ES in the parallel environments to the
performance in a sequential environment.
We use the number of {\em runs} as the $y$-axis in this figure. One run evaluates a number of potential solutions up to the number of processors. In an iteration, if the number of potential solutions is less than the number of processors, then all the potential solutions are evaluated in a single run, with some processors idle. The number of runs is equal to the number of simulations if we have only one processor.

In Fig. \ref{fig:case1_p}, we compare the number of runs needed to get from the initial NPV to 50\% of the final NPV, as well as the number of runs needed to reach the maximum number of simulator evaluations (100 times the problem dimension). From this figure we can see that, with the increase of processors, the number of runs required for GPS, PSO, and CMA-ES decrease, until the number of processors is larger than the population size. For an algorithm, the larger the population size, the greater the benefits from the parallel environment. 
With an increase in the number of processors, the order of three algorithms changes. For Case 1A, GPS performs best followed by CMA-ES and PSO in the sequential environment (number of processors equals 1). The order becomes GPS$>$PSO$>$CMA-ES when the number of processors reaches 32. Furthermore, with enough processors ($\geq 100$), the order becomes PSO$>$GPS$>$CMA-ES. The order also changes depending on the number of processors for Case 1B--1D. 
Generally, A parallel environment can greatly reduce the time spent for these algorithms. PSO can outperform GPS and CMA-ES in performance if the number of processors is large enough.

\subsection{Multiscale optimization for Model 1}
\label{sec:6_3}

In this section we address the performance of the multiscale approaches (M-GPS, M-PSO, and M-CMA-ES) for well control optimization. We use the first model as described in Section \ref{sec:5_2}. 
We stop the optimization at each scale when the average relative well rate change is less than 10\% of the distance between the upper and lower bounds on the well rates.
%We use the average relative well rate change $< 10$\% of the gap between the upper and lower bound as the stopping criterion for each scale. 
The scale will no longer be refined when the relative change in the NPV is $< 10$\% between two neighboring scales. The maximum number of simulation runs for the problem is set to 3200. As M-PSO and M-CMA-ES are stochastic algorithms, 10 trials were performed for these two algorithms to assess the overall performance.

\subsubsection{Performance of the multiscale approaches}
\label{sec:6_3_1}

As in Section \ref{sec:5_1}, we consider four configurations for each multiscale approach. 
As a first test,
the multiscale optimization process is terminated when the number of control steps reaches 8 each well for configuration I--III, and 16 each well for configuration IV. 
%This choice was made based on the results in Section \ref{sec:6_1_1} which showed that Case 1D (with 32 control steps for each well) obtains less than a 1\% increase in NPV over Case 1C (with 8 control steps for each well). The multiscale approach finds good control frequencies successfully.

\begin{figure}[htbp]
\centering 
\subfigure[GPS]{ 
    \label{fig:subfig:m_gps} %% label for first subfigure 
    \includegraphics[width=0.4\textwidth]{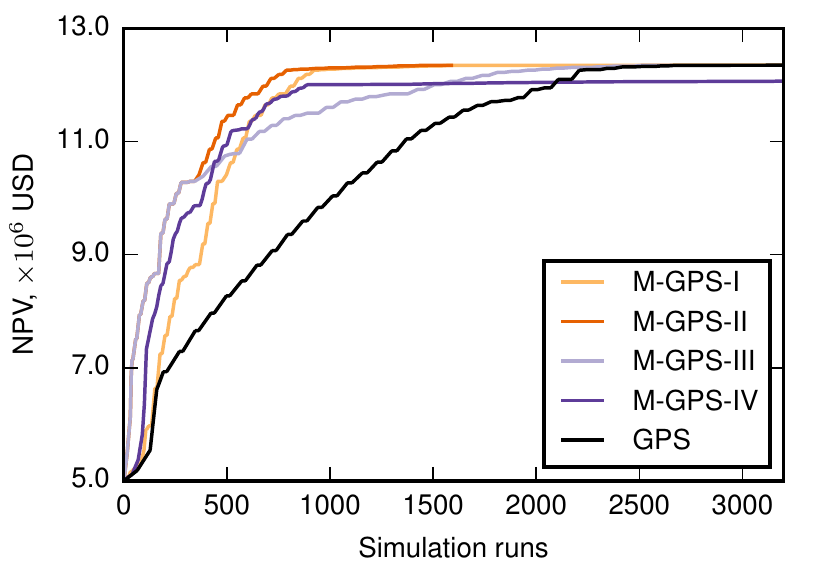}
    }
\subfigure[PSO]{ 
    \label{fig:subfig:m_pso} %% label for first subfigure 
    \includegraphics[width=0.4\textwidth]{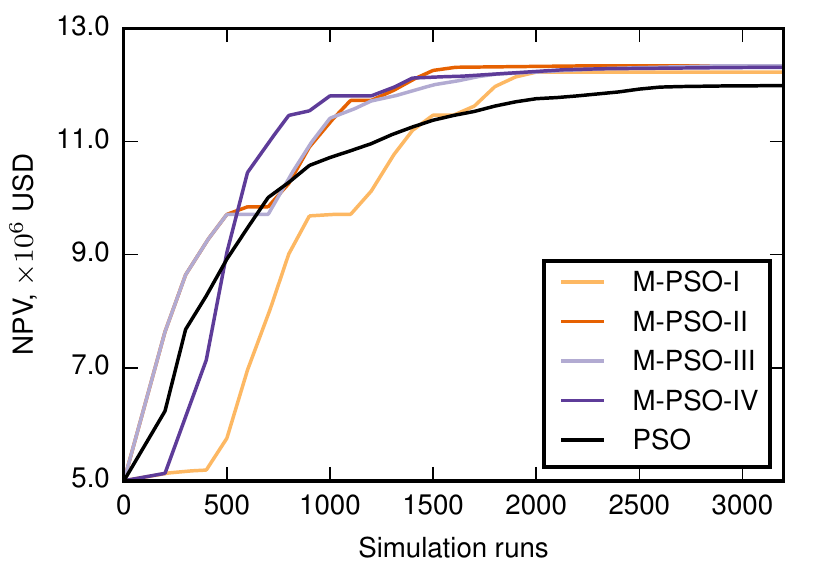}
    }
\subfigure[CMA-ES]{ 
    \label{fig:subfig:m_cmaes} %% label for first subfigure 
    \includegraphics[width=0.4\textwidth]{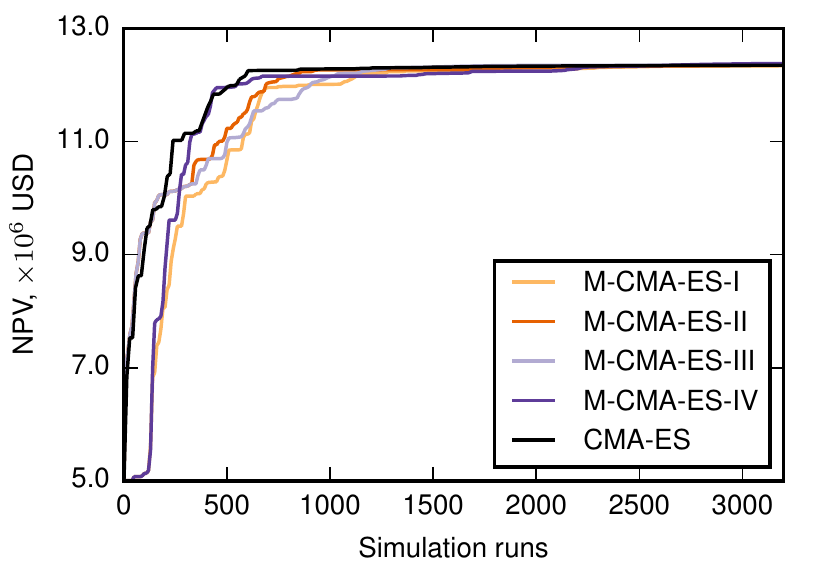}
    }
  \caption{Comparison of the performance of multiscale approaches with different configurations with 8 final control adjustments for each well.} 
  \label{fig:case1_mul} %% label for entire figure 
\end{figure}

The plots of NPV versus the number of simulation runs for the different multiscale approaches and the different configurations, as well as the plots for the standard algorithms (GPS, PSO, CMA-ES) with 8 pre-set control steps for each well, are shown in Fig. \ref{fig:case1_mul}. From the figure we see that compared with direct optimization with 8 well control adjustments, both GPS and PSO converge faster when using the multiscale approach. GPS convergence improves the most amongst the three algorithms. Fig. \ref{fig:subfig:m_pso} showed that for this test case, M-PSO could locate a control strategy which gives a higher NPV than PSO. The performance of CMA-ES (Fig. \ref{fig:subfig:m_cmaes}) is quite different. Results showed that CMA-ES converges faster than M-CMA-ES for this relatively small scale optimization problem. This is because CMA-ES is less sensitive to the quality of the initial guess and hence it can not take great advantage of the multiscale framework to speed up its convergence rate.
The multiscale framework for CMA-ES still does, however, give us a way to detect a good control frequency for well control optimization.

\begin{figure}[htbp]
\centering 
\subfigure[GPS]{ 
    \label{fig:subfig:m_gps32} %% label for first subfigure 
    \includegraphics[width=0.4\textwidth]{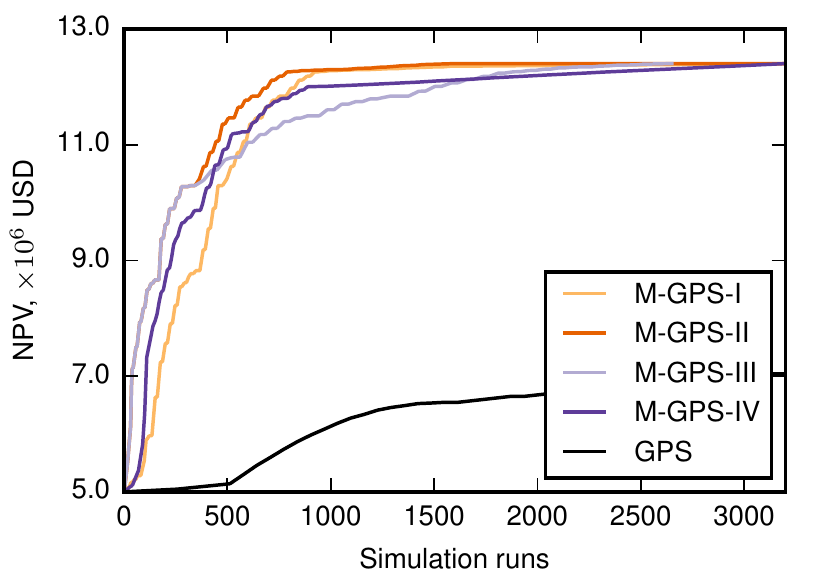}
    }
\subfigure[PSO]{ 
    \label{fig:subfig:m_pso32} %% label for first subfigure 
    \includegraphics[width=0.4\textwidth]{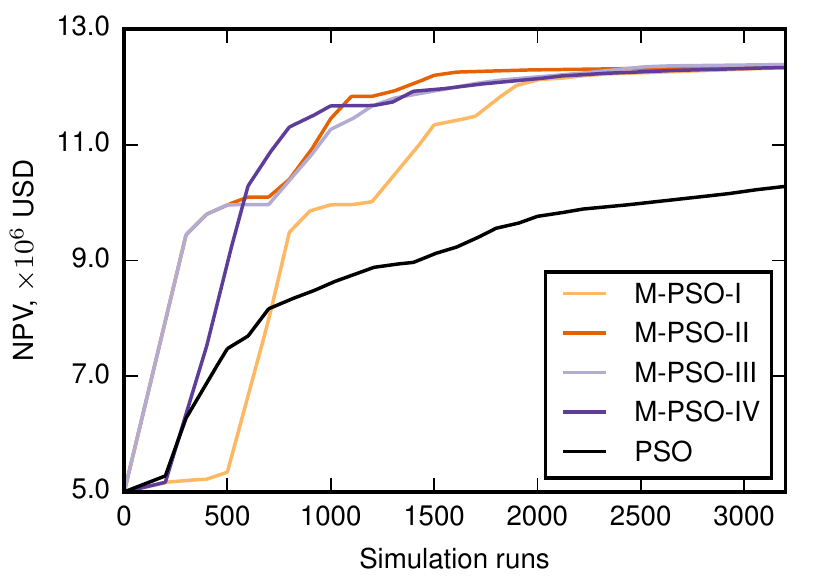}
    }
\subfigure[CMA-ES]{ 
    \label{fig:subfig:m_cmaes32} %% label for first subfigure 
    \includegraphics[width=0.4\textwidth]{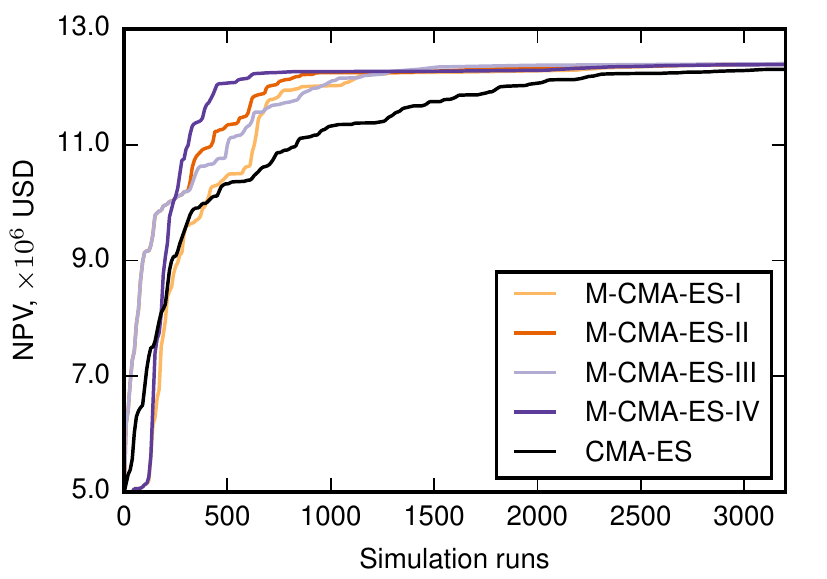}
    }
  \caption{Comparison of the performance of multiscale approaches with different configurations with 32 final control adjustments for each well.} 
  \label{fig:case1_mul32} %% label for entire figure 
\end{figure}

As a second test of the multiscale framework we increase the number of control steps to 32 for each well.
The plots of NPV versus the number of simulation runs for the different multiscale approaches and the different configurations, as well as the plots for the standard algorithms (GPS, PSO, CMA-ES) are shown in Fig. \ref{fig:case1_mul32}. From the figure we see that compared with direct optimization with 32 well control adjustments, all algorithms converge faster when using the multiscale approach. GPS convergence improves the most amongst the three algorithms, followed by PSO and CMA-ES. 

\subsubsection{Performance in parallel environments}
\label{sec:6_3_2}

\begin{figure*}[htbp]
  %\centering 
\subfigure[GPS, CPU=8]{ 
    \label{fig:subfig:m_gps_p8} %% label for first subfigure 
    \includegraphics[width=0.3\textwidth]{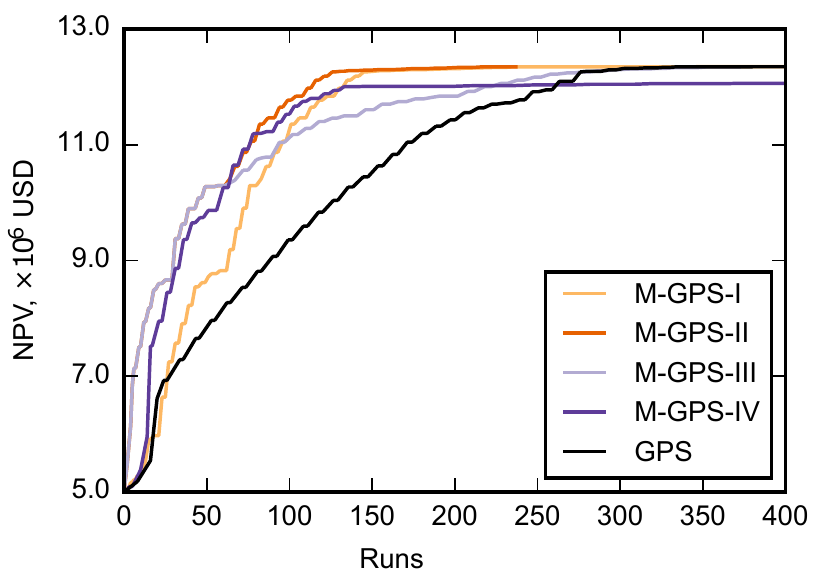}
    }
\subfigure[PSO, CPU=8]{ 
    \label{fig:subfig:m_pso_p8} %% label for first subfigure 
    \includegraphics[width=0.3\textwidth]{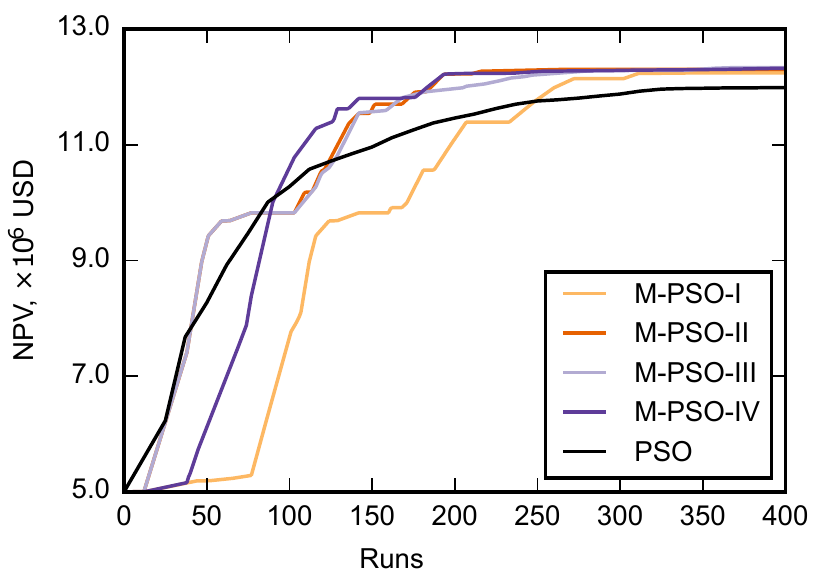}
    }
\subfigure[CMA-ES, CPU=8]{ 
    \label{fig:subfig:m_cmaes_p8} %% label for first subfigure 
    \includegraphics[width=0.3\textwidth]{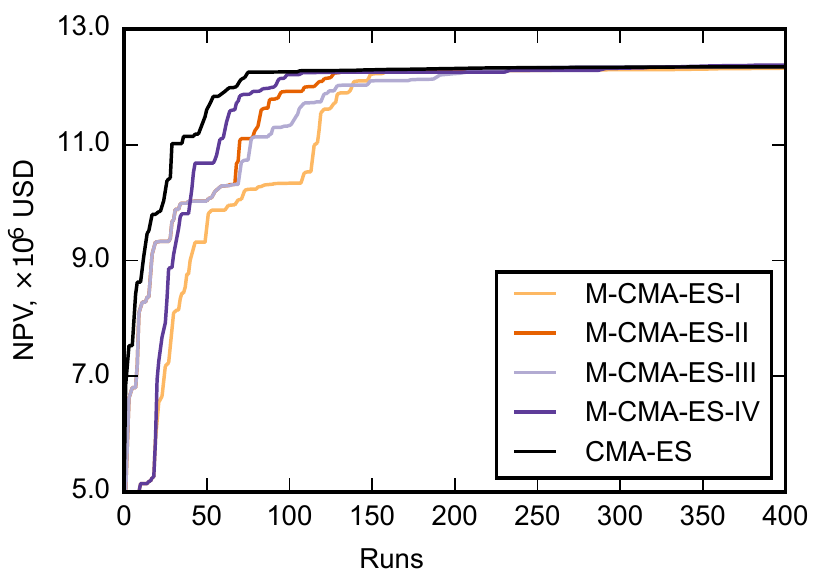}
    } 
\\
\subfigure[GPS, CPU=32]{ 
    \label{fig:subfig:m_gps_p32} %% label for first subfigure 
    \includegraphics[width=0.3\textwidth]{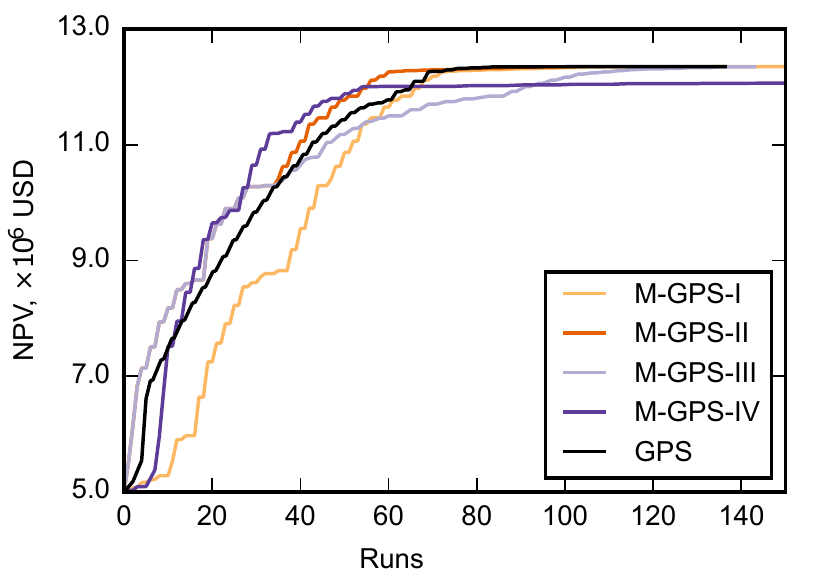}
    }
\subfigure[PSO, CPU=32]{ 
    \label{fig:subfig:m_pso_p32} %% label for first subfigure 
    \includegraphics[width=0.3\textwidth]{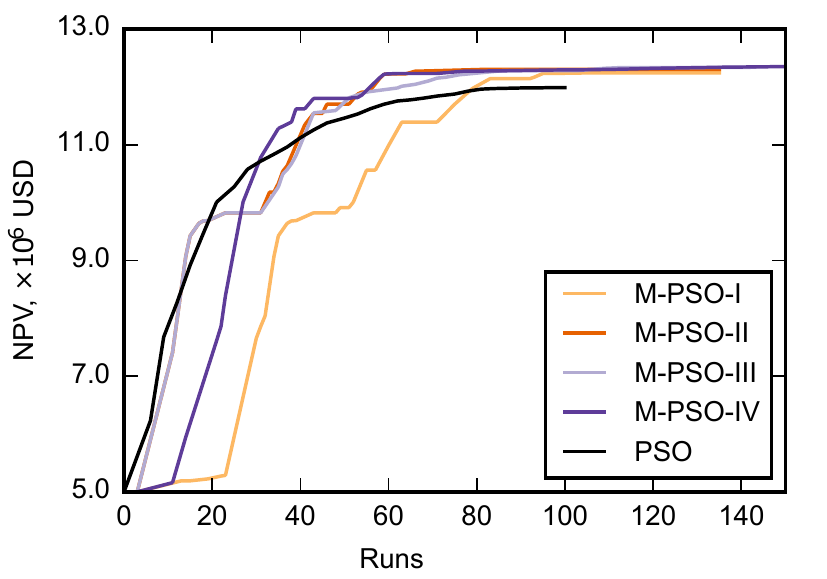}
    }
\subfigure[CMA-ES, CPU=32]{ 
    \label{fig:subfig:m_cmaes_p32} %% label for first subfigure 
    \includegraphics[width=0.3\textwidth]{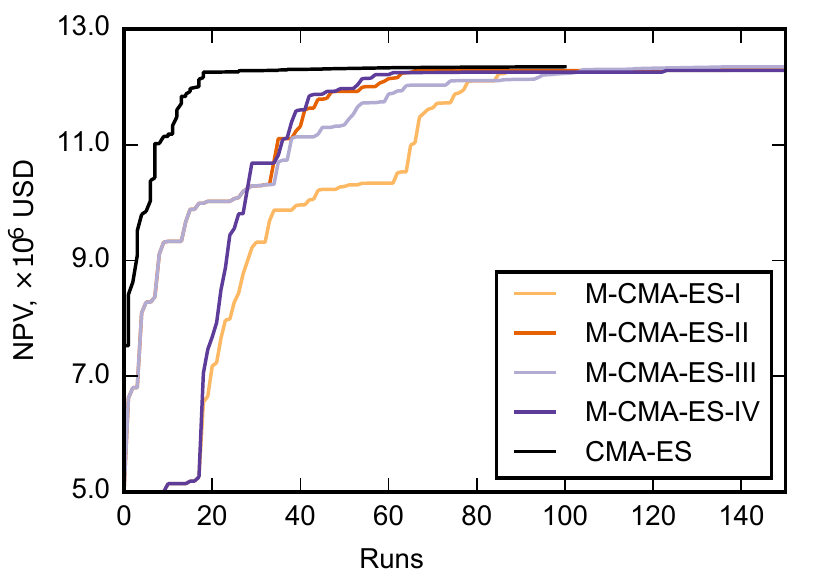}
    }
\\
\subfigure[GPS, CPU=inf]{ 
    \label{fig:subfig:m_gps_pi} %% label for first subfigure 
    \includegraphics[width=0.3\textwidth]{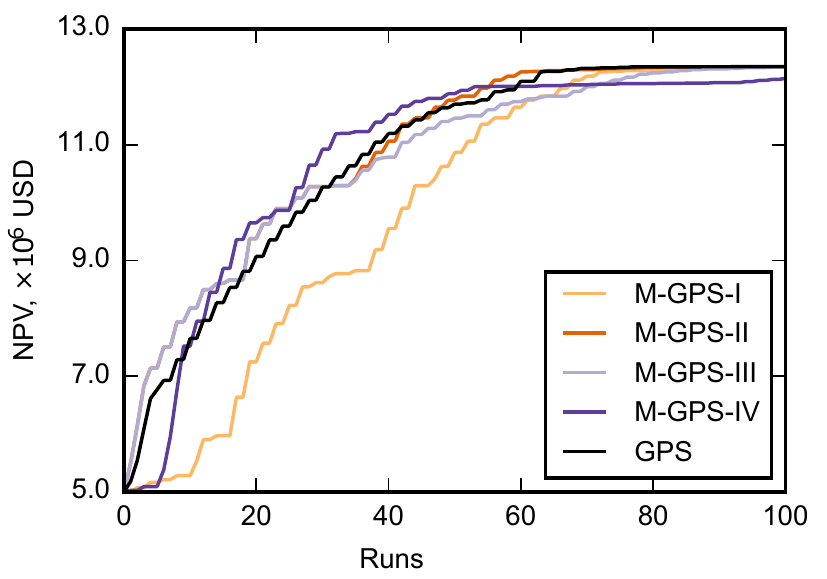}
    }
\subfigure[PSO, CPU=inf]{ 
    \label{fig:subfig:m_pso_pi} %% label for first subfigure 
    \includegraphics[width=0.3\textwidth]{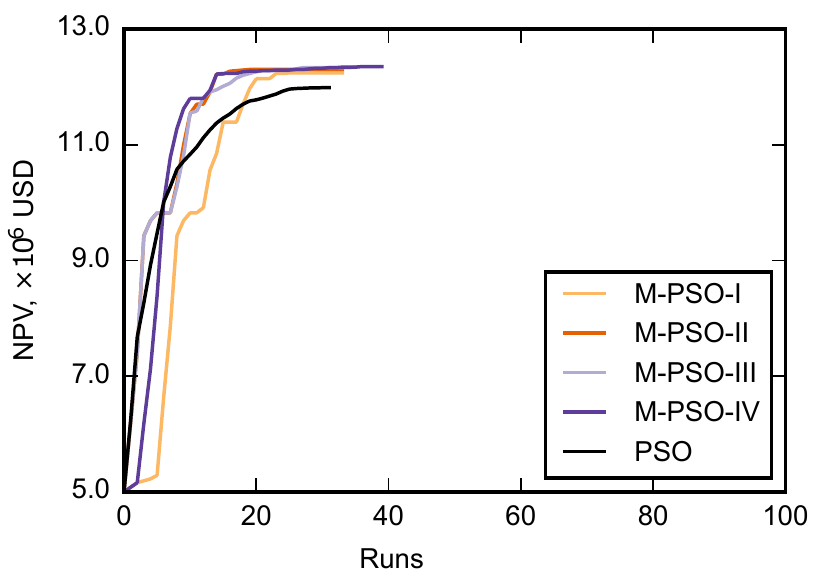}
    }
\subfigure[CMA-ES, CPU=inf]{ 
    \label{fig:subfig:m_cmaes_pi} %% label for first subfigure 
    \includegraphics[width=0.3\textwidth]{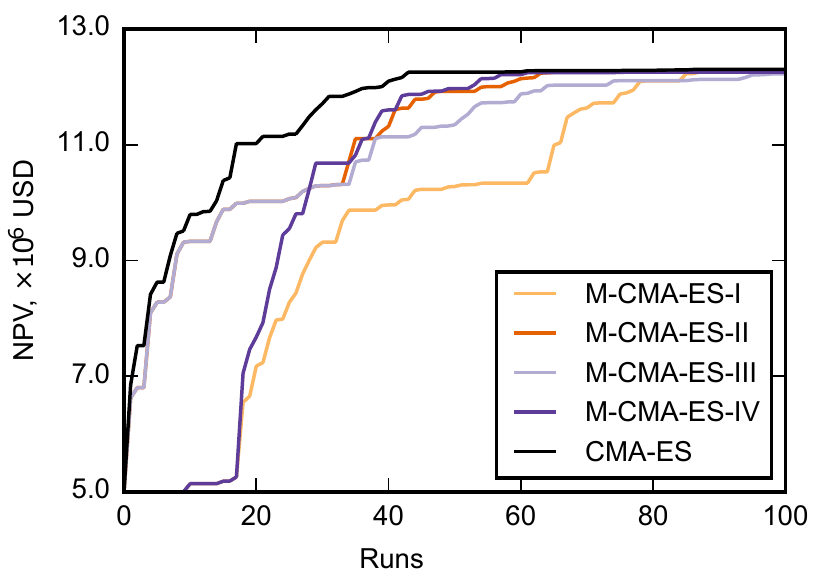}
    }
  \caption{Comparison of the performance of multiscale approaches with 8 final control adjustments for each well in parallel environments with 8, 32, and an infinite number of processors.} 
  \label{fig:case1_mul_p} %% label for entire figure 
\end{figure*}

%We investigate the performance of multiscale approaches in parallel environments. 
Fig. \ref{fig:case1_mul_p} shows the performance of the multiscale approaches in parallel environments with 8, 32, and an infinite number of processors using 8 pre-set final control steps for configurations I--III and 16 control steps for configuration IV. As in Section \ref{sec:6_2_2}, we use the number of runs as the $y$-axis. The lines for PSO, CMA-ES, M-PSO, and M-CMA-ES are the medians of 10 trials. From this figure we can see that, not surprisingly, the more processors the less runs needed to converge. In a parallel environment, the improved convergence of the multiscale approach is less apparent.
The multiscale approach still benefits, however, if the optimal control frequency is unknown and hence not specified a priori.

\begin{figure*}[htbp]
  %\centering 
\subfigure[GPS, CPU=8]{ 
    \label{fig:subfig:m_gps_p832} %% label for first subfigure 
    \includegraphics[width=0.3\textwidth]{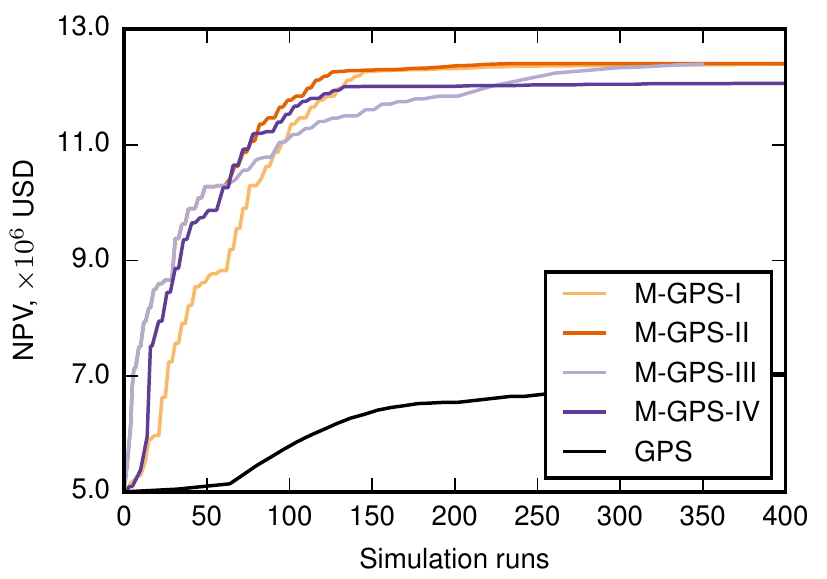}
    }
\subfigure[PSO, CPU=8]{ 
    \label{fig:subfig:m_pso_p832} %% label for first subfigure 
    \includegraphics[width=0.3\textwidth]{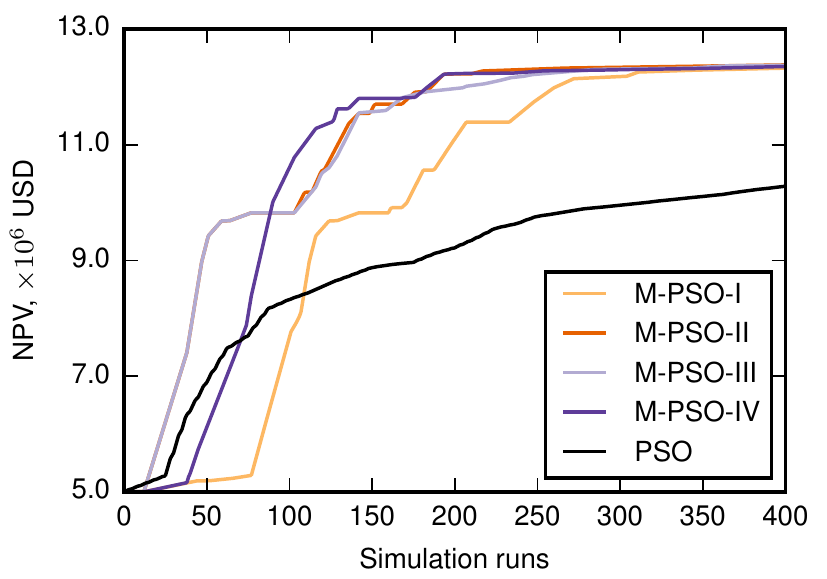}
    }
\subfigure[CMA-ES, CPU=8]{ 
    \label{fig:subfig:m_cmaes_p832} %% label for first subfigure 
    \includegraphics[width=0.3\textwidth]{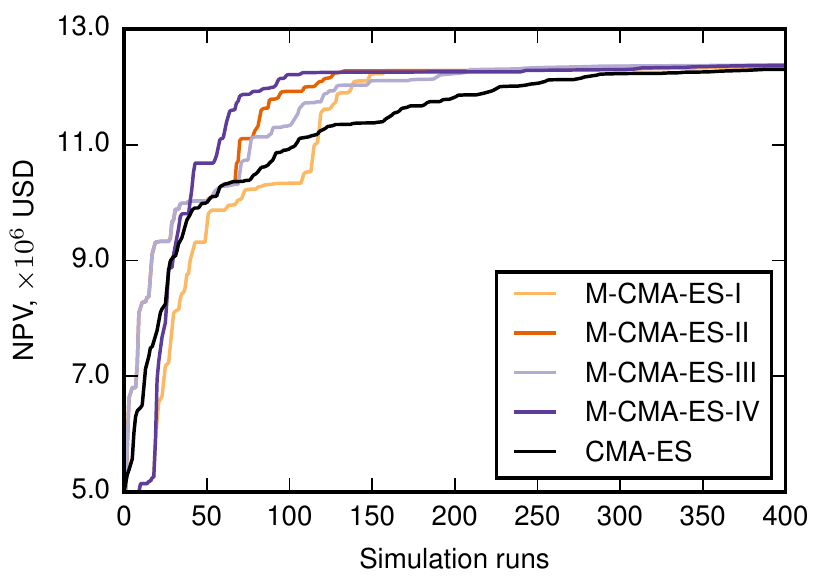}
    } 
\\
\subfigure[GPS, CPU=32]{ 
    \label{fig:subfig:m_gps_p3232} %% label for first subfigure 
    \includegraphics[width=0.3\textwidth]{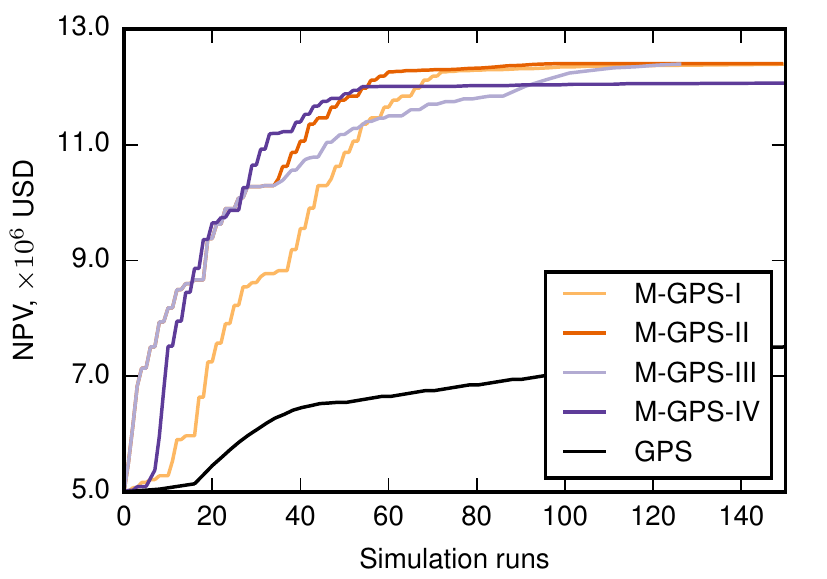}
    }
\subfigure[PSO, CPU=32]{ 
    \label{fig:subfig:m_pso_p3232} %% label for first subfigure 
    \includegraphics[width=0.3\textwidth]{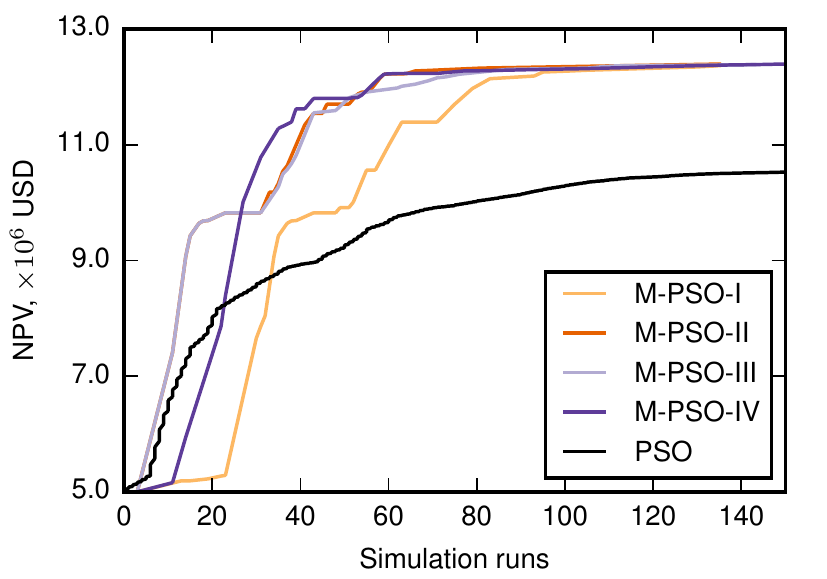}
    }
\subfigure[CMA-ES, CPU=32]{ 
    \label{fig:subfig:m_cmaes_p3232} %% label for first subfigure 
    \includegraphics[width=0.3\textwidth]{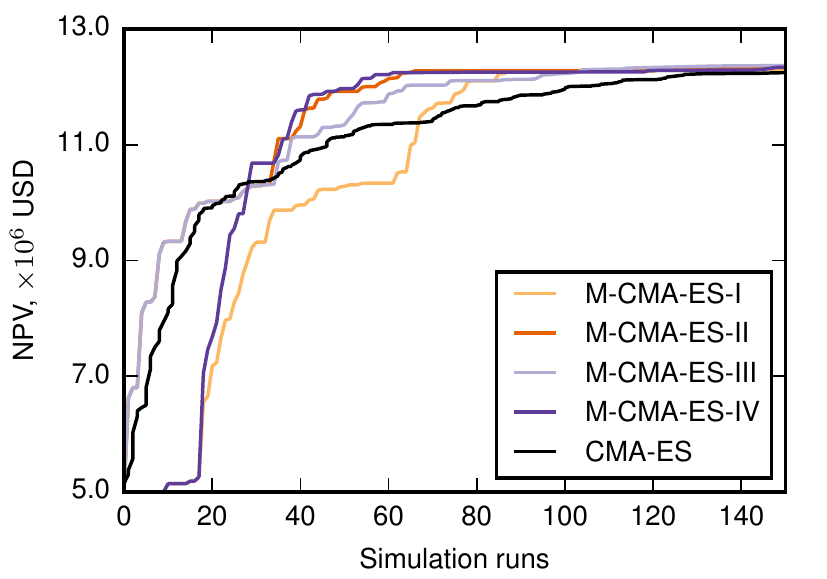}
    }
\\
\subfigure[GPS, CPU=inf]{ 
    \label{fig:subfig:m_gps_pi32} %% label for first subfigure 
    \includegraphics[width=0.3\textwidth]{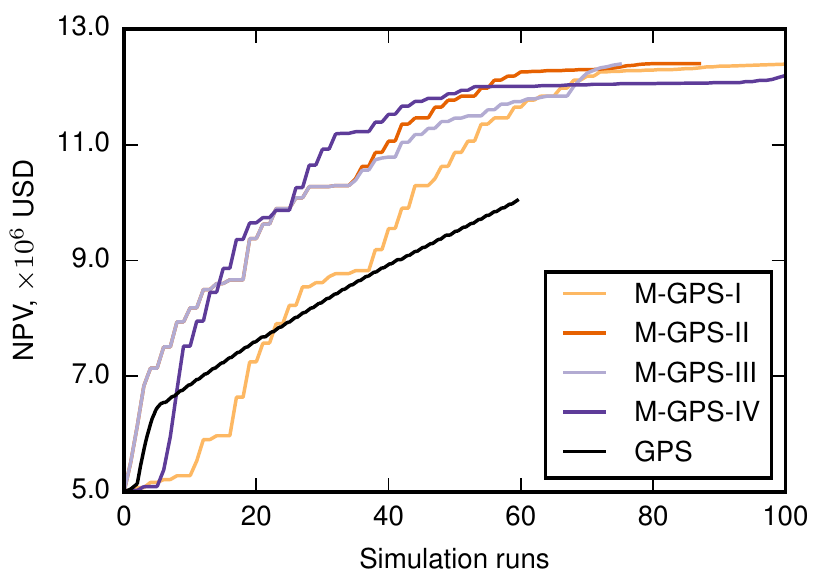}
    }
\subfigure[PSO, CPU=inf]{ 
    \label{fig:subfig:m_pso_pi32} %% label for first subfigure 
    \includegraphics[width=0.3\textwidth]{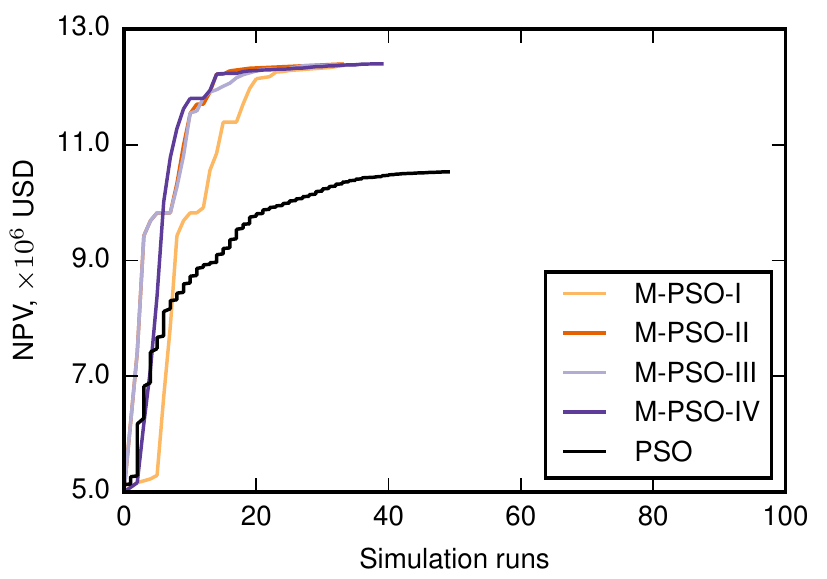}
    }
\subfigure[CMA-ES, CPU=inf]{ 
    \label{fig:subfig:m_cmaes_pi32} %% label for first subfigure 
    \includegraphics[width=0.3\textwidth]{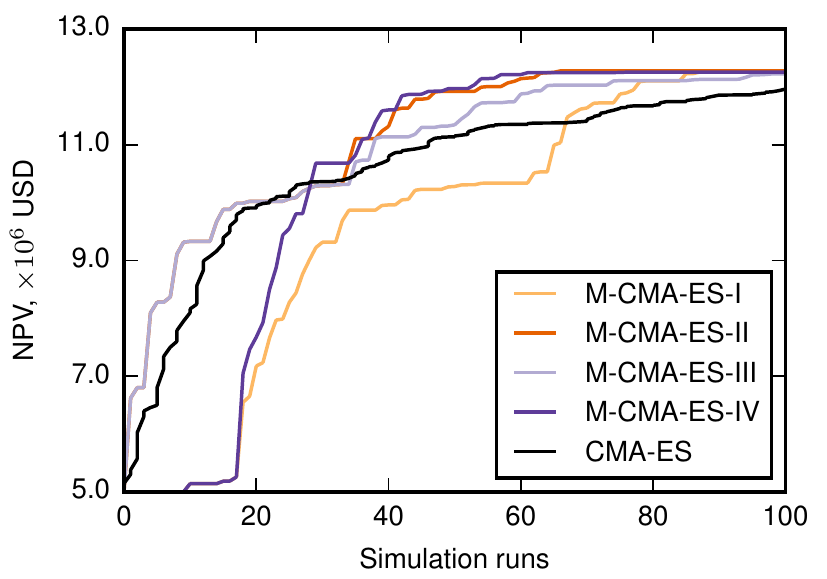}
    }
  \caption{Comparison of the performance of multiscale approaches with 32 final control adjustments for each well in parallel environments with 8, 32, and an infinite number of processors.} 
  \label{fig:case1_mul_p32} %% label for entire figure 
\end{figure*}

In Fig. \ref{fig:case1_mul_p32} we repeat this experiment for the more difficult problem with 32 control steps for each well. The results are similar.

\subsubsection{Performance with different computational budgets}
\label{sec:6_3_3}

\begin{figure*}[htbp]
\centering 
\subfigure[Low budget, 300 simulation runs]{ 
    \label{fig:subfig:m_low} %% label for first subfigure 
    \includegraphics[width=.8\textwidth]{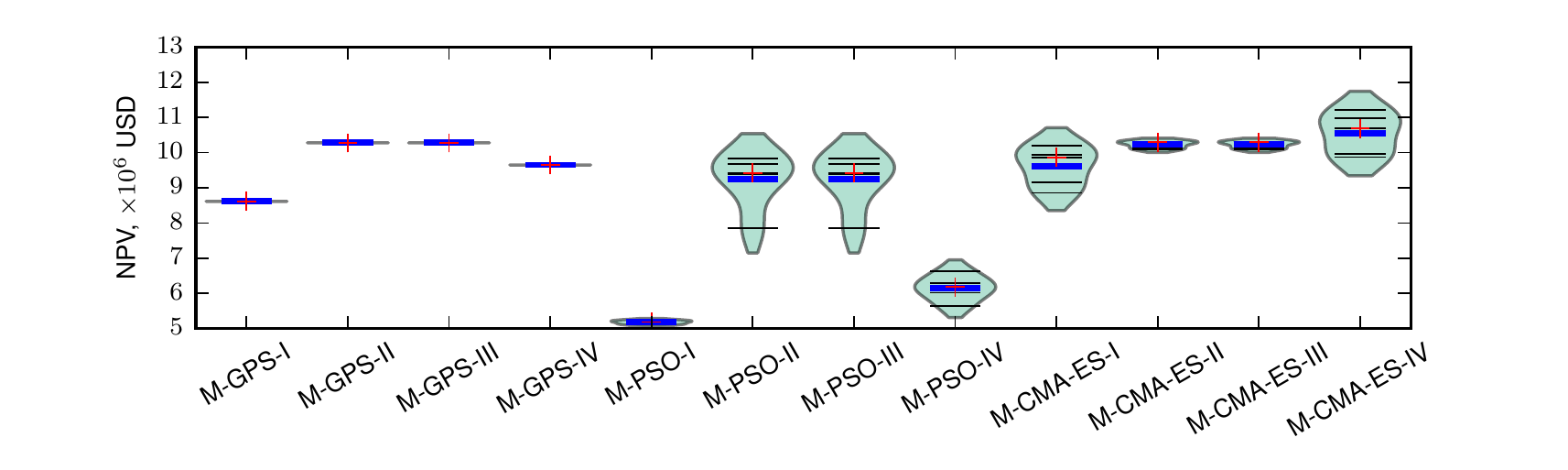}
    }
\subfigure[Medium budget, 1500 simulation runs]{ 
    \label{fig:subfig:m_mid} %% label for first subfigure 
    \includegraphics[width=.8\textwidth]{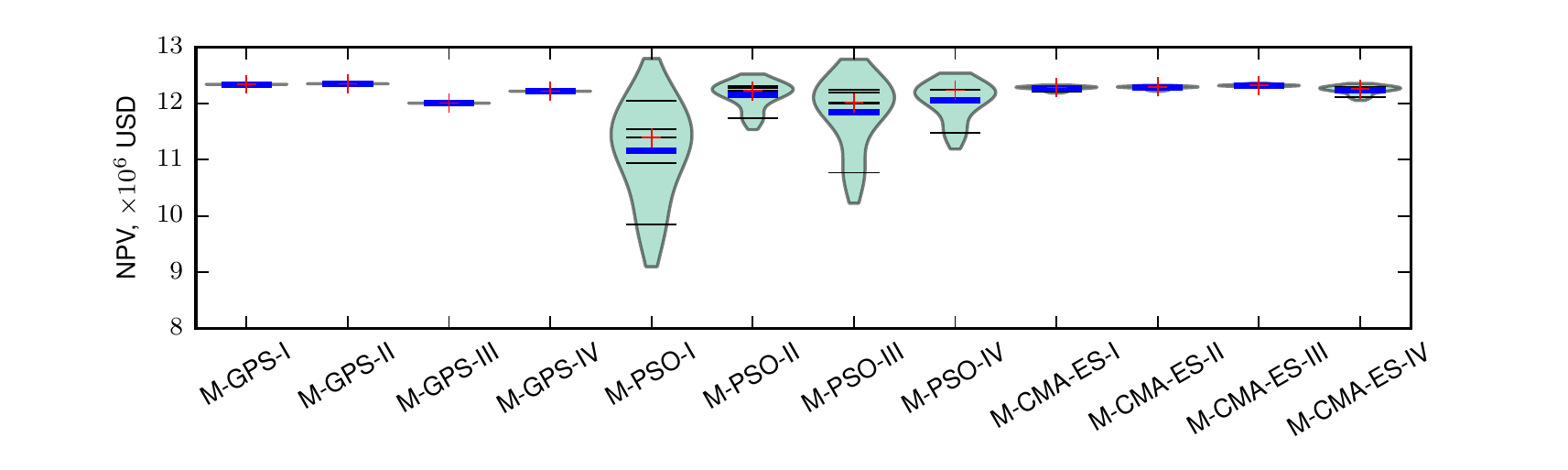}
    }
\subfigure[High budget, 3000 simulation runs]{ 
    \label{fig:subfig:m_high} %% label for first subfigure 
    \includegraphics[width=.8\textwidth]{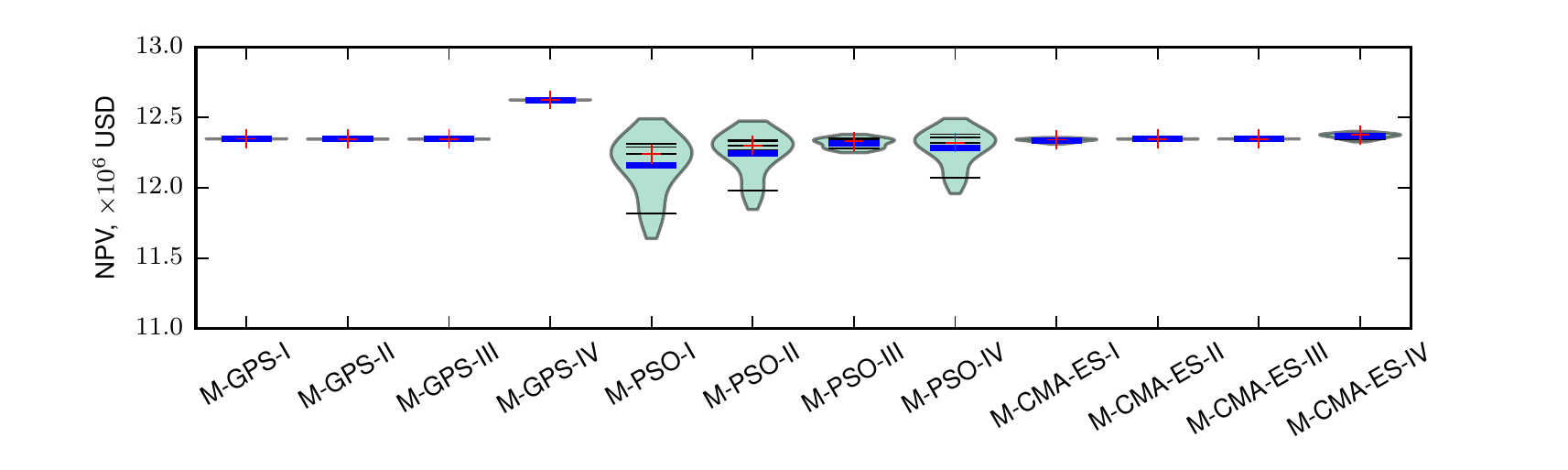}
    }
  \caption{Beanplots for different configurations of our multiscale approaches.} 
  \label{fig:e2_beanplot} %% label for entire figure 
\end{figure*}

Here we assess the effect of computational budget on the efficacy of the multiscale approach. We use 300, 1500, and 3000 simulation runs as a low, medium, and high budget. For different budgets, the performance of the different multiscale approaches and the different configurations are shown in several beanplots in Fig. \ref{fig:e2_beanplot}. In the beanplots, the individual observations are shown as small lines in a one-dimensional scatter plot. The estimated density of the distributions is visible and the mean (blue bold line) and median (red `+') are shown.

Fig \ref{fig:subfig:m_low} shows results for a low budget. Since GPS is a deterministic algorithm, we have no distribution for the four configurations of M-GPS. M-GPS-II and M-GPS-III obtained the highest NPV amongst all four configurations. Although some configurations of M-PSO and M-CMA-ES could obtain a relatively high NPV, there is also has a risk of obtaining a low NPV due to the high variability. For a medium budget, Fig \ref{fig:subfig:m_mid}, the variation of all four configurations of M-CMA-ES is relatively small. The variation of the M-PSO is quite large. With this budget of function evaluations, M-GPS and M-CMA-ES are good choices. For a large budget, we see that M-GPS-IV obtained a higher NPV than the other configurations. With the fourth configuration, the approach terminated at 16 control steps for each well, while the other three configurations terminated at 8 control steps for each well. With more control steps, a higher NPV could be obtained. Configuration I performs less well with a low budget. This is because the initial number of control steps for configuration I is 1. The optimal control found with only 1 step is very different when compared with the optimal control found with more steps. In general, we see that the second configuration performs better than the other three configurations for M-GPS, M-PSO and M-CMA-ES. Furthermore, M-GPS-II is highly recommended for all budgets.   

\subsection{Multiscale optimization for a real-world reservoir}
\label{sec:6_4}

Based on the results of the previous section we now apply the multiscale approaches with configuration II ($n_0=2$ and $n_s=2$) to solve the well control problem of reservoir PUNQ-S3 (see Section \ref{sec:5_2_2}). The maximum number of function evaluations is set to 10000. We use an average relative well rate change of $< 10$\% of the gap between the upper and lower bound as the stopping criterion for each scale. The scale will no longer be refined when the relative change in the NPV is $< 10$\% between two neighboring scales. Due to the large computational time we perform only $3$ trials for M-PSO and M-CMA-ES.

\begin{figure}[htbp]
\centering 
\subfigure[GPS]{ 
    \label{fig:subfig:2m_gps} %% label for first subfigure 
    \includegraphics[width=0.4\textwidth]{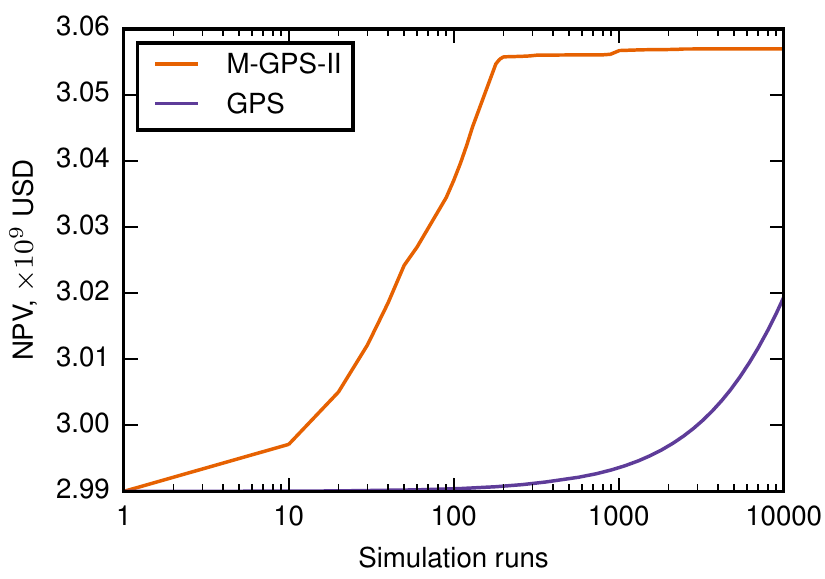}
    }
\subfigure[PSO]{ 
    \label{fig:subfig:2m_pso} %% label for first subfigure 
    \includegraphics[width=0.4\textwidth]{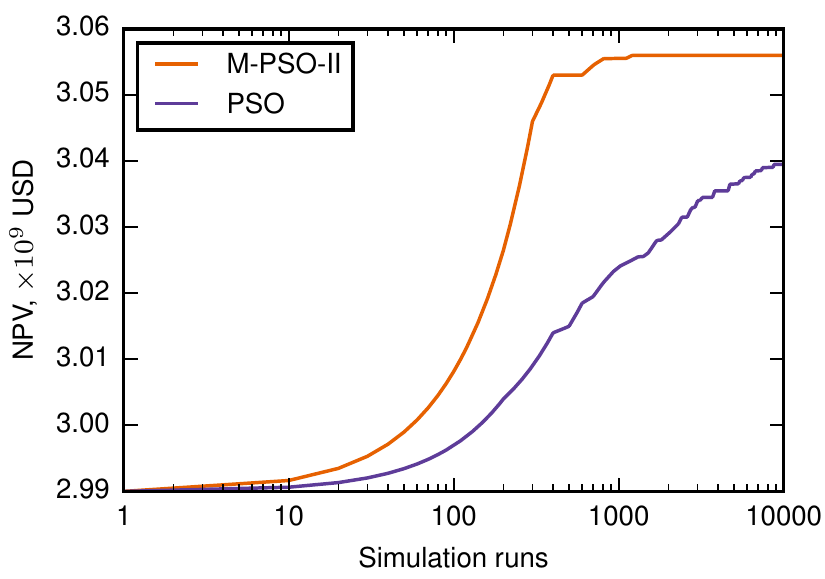}
    }
\subfigure[CMA-ES]{ 
    \label{fig:subfig:2m_cmaes} %% label for first subfigure 
    \includegraphics[width=0.4\textwidth]{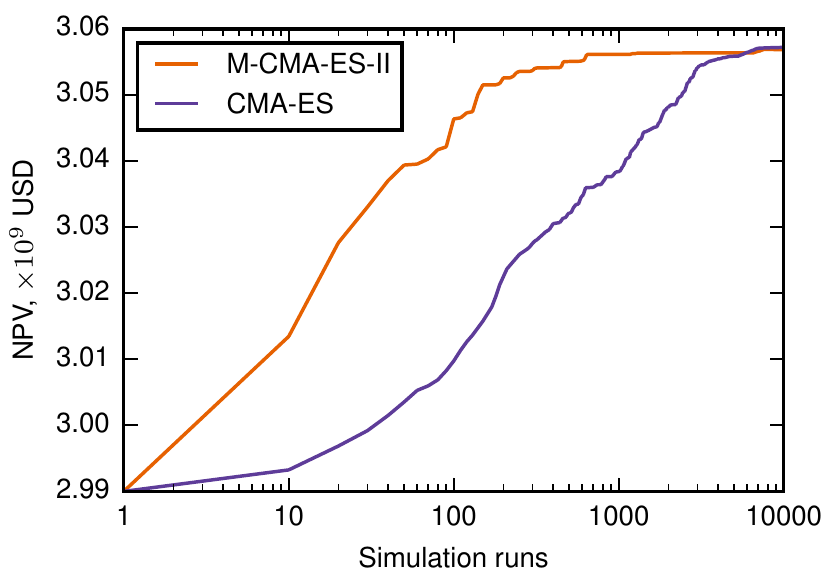}
    }
  \caption{Comparison of the performance of multiscale approach and optimizers without multiscale for PUNQ-S3. Here are the median NPV of trials for PSO and CMA-ES. } 
  \label{fig:punqs3_mul} %% label for entire figure 
\end{figure}

Fig. \ref{fig:punqs3_mul} shows the performance of M-GPS-II, M-PSO-II and M-CMA-ES-II. The performance of GPS, PSO, and CMA-ES using a pre-set control frequency of 32 control steps for each well is shown in this figure as well. The results show that for the same number of reservoir simulations, combining these three algorithms with the multiscale framework gives higher NPV values as compared to directly optimizing with the largest number of control steps. 

The values of NPV found by the approaches after 1000, 5000, and 10000 simulation runs are shown in Table \ref{tab:punqs3_result_npv}. From the table we can see that, M-GPS-II outperforms among all three multiscale approaches, even when the number of simulation runs is limited. Of course the deterministic M-GPS has the additional advantage of having no variability in the outcome. Without the multiscale framework, CMA-ES performs the best. 

\begin{table}[htbp]
\centering 
\caption{Median NPV ($10^9$ USD) with different budgets for PUNQ-S3.}
\label{tab:punqs3_result_npv}
\begin{tabular}{lcccc}\hline
Algorithm   & Trails & 1000   & 5000   & 10000  \\\hline
GPS         & 1     & 2.9936 & 3.0061 & 3.0194 \\
M-GPS-II    & 1     & 3.0567 & 3.0570 & 3.0570 \\\hline
PSO         & 3     & 3.0240 & 3.0365 & 3.0395 \\
M-PSO-II    & 3     & 3.0555 & 3.0560 & 3.0560 \\\hline
CMA-ES      & 3     & 3.0384 & 3.0558 & 3.0572 \\
M-CMA-ES-II & 3     & 3.0561 & 3.0563 & 3.0568 \\\hline
\end{tabular}
\end{table}

\section{Conclusions}
\label{sec:7}

In this paper we have considered three derivative-free optimization algorithms combined with a multiscale framework to solve well control optimization problems. The optimization algorithms used include GPS, which is a deterministic local search approach; PSO, which is a stochastic global search method; and CMA-ES, which is a stochastic local search method. A generalization of the successive-splitting multiscale approach from \cite{lien_multiscale_2008,shuai_using_2011} was introduced to combine with the derivative-free optimization algorithms. 

Based on thorough numerical experiments the following conclusions can be drawn:
\begin{itemize}
\item[\textbullet] The control frequency does have a significant effect on well control optimization problems. The more frequent the well control adjustment, the higher the NPV that can be obtained but at the cost of a harder optimization problem.
This increase becomes less significant as we continue to increase the number of control steps.
%but the increase of NPV becomes less. 
Considering the operation costs, each reservoir has a optimal control frequency. The optimal controls are similar with different control frequencies when every well is produced under a liquid rate throughout its lifetime. 
\item[\textbullet] Without the multiscale framework,
%For the performance of GPS, PSO, and CMA-ES (without combining the multiscale approach) in solving well control problem, 
GPS performs best when the problem dimension is very small and the budget is large enough. CMA-ES showed excellent performance when the budget is limited. A parallel environment can greatly reduce the time spent for these algorithms. PSO can outperform GPS and CMA-ES in performance if the number of processors is large enough. The choice of the initial guess has a significant effect on the convergence speed for GPS and PSO at the early stage of optimization. CMA-ES, by contrast, is less sensitive to the choice of initial guess. 
%For PSO, another parameter that influences the performance is the population size.
The performance of PSO is affected dramatically by the population size.
Parameter tuning for CMA-ES showed that the default settings work quite well.
\item[\textbullet] The multiscale approaches have two advantages in solving well control problem. (1) they provide a way to optimize the control frequency and the well controls simultaneously. (2) when compared to the standalone algorithms the multiscale approach can speed-up the convergence. Based on the results of the test cases, the convergence of GPS and PSO improves the most when combined with the multiscale framework.
The multiscale framework is more efficient as the number of control steps increases.
%The improved speed-up of the multiscale approach is less obvious in a parallel environment.
%The advantage that speed-up the convergence rate for the multiscale approaches become less obvious in the parallel environment.
The difference in performance between the multiscale hybrid algorithms and the stand-alone algorithms decreases as the number of processors increases.
\item[\textbullet] The multiscale framework has two key parameters, the choice of the initial number of control steps $n_0$ and the split factor $n_s$. The choice $n_0=2$ and $n_s=2$ gave the best performance. In the multiscale framework, M-GPS-II is highly recommended for any computational budget.
\end{itemize}

All above conclusions are based on the experiments in this paper. Although the multiscale approach has shown its potential to solve complex well control problems, there are still many possible avenues for future work. Some potential areas of investigation are the use of other stochastic approaches such as mesh adaptive direct search or differential evolution to see if they offer any improvement in performance. There is also flexibility to choose different $n_0$ and $n_s$ values for each well. The performance of multiscale approaches with nonlinear constraints still needs additional study. As does the development of robust stopping criteria within the multiscale framework.

\begin{acknowledgements}
The authors acknowledge funding from the Natural Sciences and Engineering Research Council of Canada (NSERC) Discovery Grant Program, the National Science and Technology Major Project of the Ministry of Science and Technology of China (2011ZX05011-002), and the program of China Scholarships Council (No. 201406450017).
%the Program for Changjiang Scholars and  Innovative Reserch Team in University (IRT1294).
\end{acknowledgements}

%\bibliography{Refexpand}   % name your BibTeX data base

% % Non-BibTeX users please use
% \begin{thebibliography}{}
% %
% % and use \bibitem to create references. Consult the Instructions
% % for authors for reference list style.
% % Format for Journal Reference
% Author, Article title, Journal, Volume, page numbers (year)
% Author, Book title, page numbers. Publisher, place (year)
% \end{thebibliography}

\end{document}